\magnification\magstep1
\input epsf.sty
\def\cC{{\cal C}}
\def\cI{{\cal I}}
\def\cX{{\cal X}}
\def\C{{\bf C}}
\def\P{{\bf P}}
\def\R{{\bf R}}
\def\Z{{\bf Z}}
\font\contentsrm=cmr10 at 8pt

\centerline{\bf Dynamics of a Two Parameter Family }

\centerline{{\bf of Plane  Birational Maps:  Maximal entropy}\footnote*{The authors thank the NSF for support during the preparation of this paper.}}
\medskip
\centerline{Eric Bedford and Jeff Diller}
\medskip
{\obeylines\contentsrm
Contents:
0.  Introduction
1.  Family of Examples: Regularization of the Generic Map
2.  Action of the Pullback on Divisors
3.  Real Mappings and their Compactifications
4.  Combinatorics of Real Curves
5.  Coding Orbits
6.  Orbits Attracted to Infinity
7.  Invariant Measure}

\bigskip\noindent {\bf \S0.  Introduction. }
Here we study a family of plane birational mappings $f_{a,b}$ which was introduced in [BMR] in connection with infinite discrete symmetries in two-dimensional lattice statistical mechanics models.  This family has been studied further
by several authors, including the works [BHM1] and  [A1--3].   These papers identify a number of properties, largely obtained from numerical investigations, and raise a number of questions.  First among these (cf [BV]) is:  How do the degrees of the iterates $f^{n}_{a,b}=f\circ\cdots\circ f$ grow as $n\to\infty$?  They found that for generic values of $a$ and $b$, the degree growth is asymptotically $\rho^{n}$ with $\rho\sim 2.1479$.  They also observed some values of $a$ and $b$ for which the growth rate appeared to be smaller.  When the parameters $a$ and $b$ are both real, $f_{a,b}$
defines a birational map of ${\bf R}^{2}$.  The second line of inquiry  pursued by these authors was driven by a desire to understand aspects of the dynamics of the family of real mappings $f_{a,b}$.
\bigskip
\epsfysize=2in
\centerline{ \epsfbox{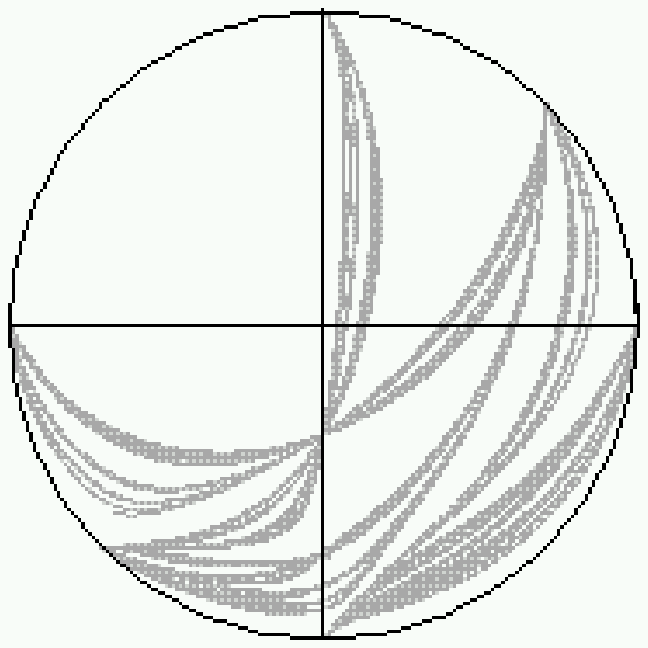} \hfil\epsfysize=2in\epsfbox{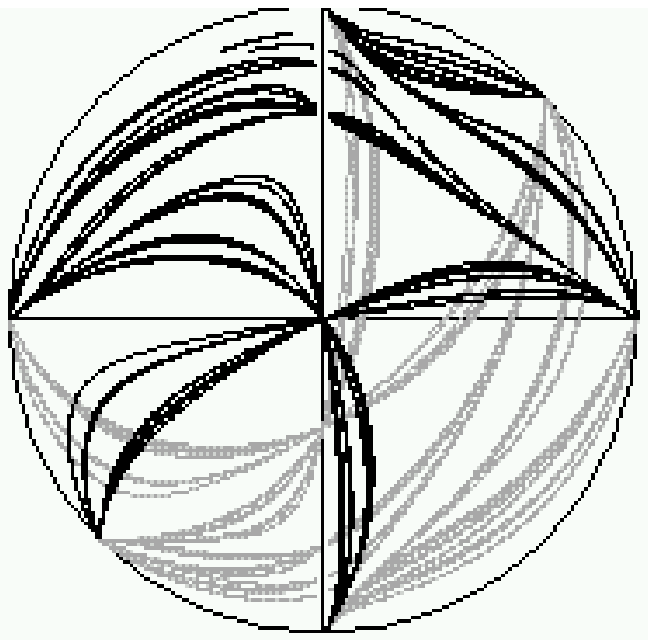}  }
\centerline{Figure 0.1.  $f_{a,b}$ with $a=-2$, $b=1$; an unstable manifold (left),}

\centerline{stable and unstable manifolds (right).}
\medskip

This paper re-examines these mappings and these questions.  Earlier, in [BD1], we made a related study.  That paper involved a family which is birationally equivalent to the the sub-family $f_{a,0}$, corresponding to $b=0$.  The degree growth of that family for generic $a$ is given by the golden mean $\phi\sim1.61$.   Each mapping $f_{a,0}$ has a parabolic fixed point in the non-wandering set, and is thus not hyperbolic.   In [BD1] we analyzed mappings $f_{a,0}$ with (maximal) entropy equal to $\log\phi$.  We showed how, when the entropy is maximal, complex methods could be used in lieu of hyperbolicity.  We showed (1)  that the non-wandering set is the complement of the parabolic basin, and (2) that the restriction to the non-wandering set was essentially conjugate to the golden mean sub-shift.

In this paper we address the general case $b\ne0$.   To illustrate one of these mappings, we have drawn in Figure 0.1 the unstable manifold $W^{s}(p)$ of a saddle point $p$ for $f_{a,b}$ with $a=-2$, $b=1$.  This is presented in coordinates $(\rho,\theta)$, which are obtained from the usual polar coordinates $(r,\theta)$ by setting $\rho={2\over\pi}\arctan r$.  Thus the disk  $\{\rho<1\}$ corresponds to ${\bf R}^{2}$, and $\{\rho=1\}$ is the circle at infinity.  The intersection $W^{s}(p)\cap W^{u}(p)$, as seen in the right half of this figure, suggests a large set of homoclinic points and leads us to expect that $f$ will have interesting behavior.  One difference between the cases $b\ne0$ and $b=0$ may be seen by contrasting the pinching at infinity in Figures 0.1 and 0.2.  These two cases also exhibit different expansion/entropy, which is reflected in the lengths of the curves, since in both cases what is drawn is $f^{9}I$, where $I$ is an arc which crosses the disk once.
\bigskip
\epsfysize=2in
\centerline{ \epsfbox{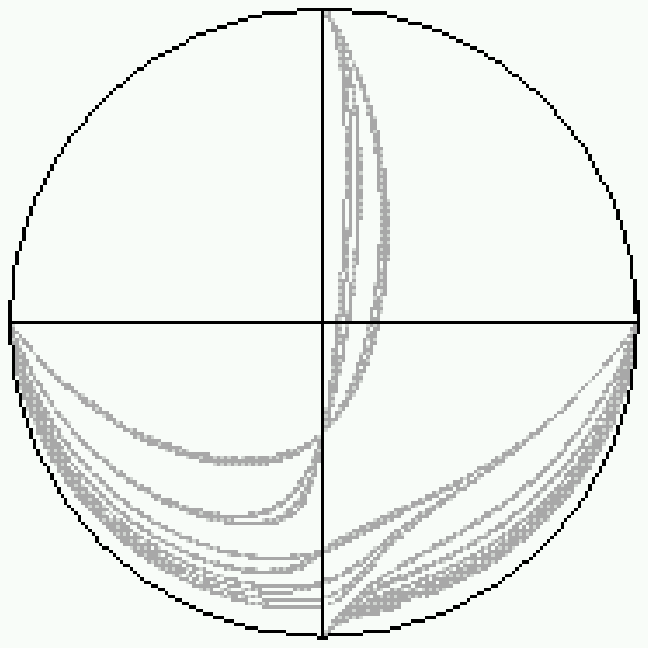} \hfil\epsfysize=2in  \epsfbox{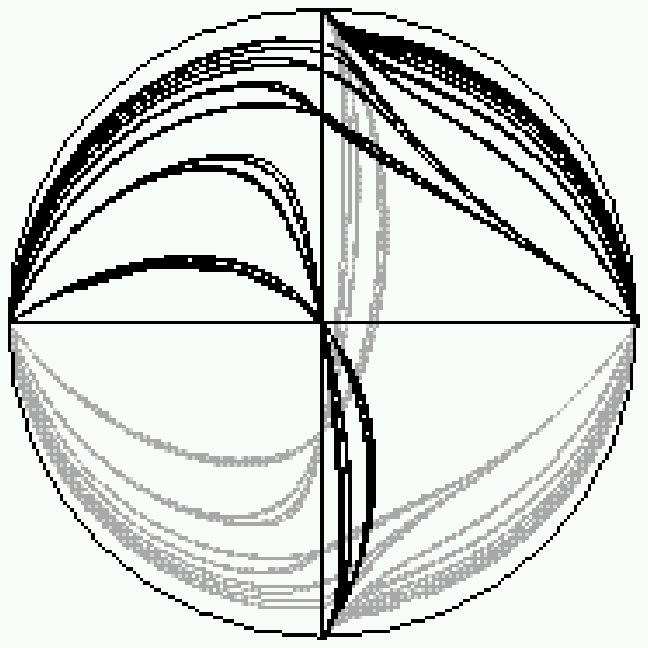} }
\centerline{Figure 0.2.  $f_{a,b}$ with $a=-2$, $b=0$; an unstable manifold (left),}

\centerline{stable and unstable manifolds (right).}
\medskip

A central strategy in this paper is to look at the action of $f_{a,b}$ on one-dimensional objects, both real and complex.  On the complex side, one associates to any compact complex manifold $\cX$ its Picard group ${\rm Pic}(\cX)$, which is the set of divisors on $\cX$ modulo linear equivalence.  A rational map between compact complex surfaces induces linear pullback and pushforward actions $f^*,f_*$ on ${\rm Pic}(\cX)$. In \S1 our primary concern is to construct a complex surface $\cX$ on which the property $(f^*)^n = (f^n)^*$ holds.  To do this we follow the method of [BTR] and [DF] and blow up the orbits of exceptional curves.

In \S2, we compute the action $f^{*}$ on ${\rm Pic}(\cX)$, and we verify (Theorem 2.1) that for generic parameters
$a,b\in{\bf C}$ the degree growth rate is exponential and given by  $\rho\sim2.1479$, which is the largest root of the
polynomial $x^{3}-x^{2}-2x-1$. When the parameters $a$ and $b$ are real, the set $\cX_R := \overline{\R^2}\subset \cX$ is an $f$-invariant compact real 2-manifold.  In \S3 we discuss the topology of $\cX_\R$.  This leads us to a useful combinatorial device (the ``dynamic hexagon'' in Figure 3.3) for analyzing the real dynamics of $f$.

In \S4 we identify two classes of real arcs in $\cX_\R$ on which the actions of $f$ and $f^{-1}$ may be analyzed combinatorially.  Using combinatorial analysis and complex intersection theory, we show that these real actions are essentially equivalent to the complex actions $f^*$ and $f_*$.  Pushing our analysis further in \S5, we define a `coding map' $c:\Omega\to \Sigma_F$ from an $f$-invariant set $\Omega\subset\cX_\R$ onto a subshift $\Sigma_F$ of finite type.  In \S6 we show that the complement of $\Omega$ in $\cX_\R$ is the basin of infinity.

One of our goals is to make the connection between the real and complex dynamics of $f$.  Our conclusion is that when $a<-1$, $b\ne0$, the real and complex dynamics are essentially the same.  We show in Theorem 6.4 that all homoclinic points are real: if $W^s_C$ and $W^u_C$ are complex stable manifolds of a real saddle point, then $W^s_C\cap W^u_C\subset{\bf R}^2$.  In \S7 we apply some results from the theory of birational dynamics of compact complex surfaces.  In particular, the general theory gives the existence of a canonical invariant measure $\mu$.  In Theorem 7.1, we show that $\mu$ is carried by $\Omega\subset{\bf R}^{2}$.  Then, we show (Theorem 7.2) that the pushforward $c_*\mu$ is the unique measure of maximal entropy on $\Sigma_F$.

\bigskip\noindent{\bf \S1.  Family of Examples: Regularization of the Generic Map. } The mappings we consider are of the form $f=\tau\circ\sigma$, where
$$
\eqalign{ \sigma(x,y)=(1-x+{x\over y},1-y+{y\over x}),\cr
\tau(x,y)=(x,b x + a + 1-y).\cr} \eqno(1.1)
$$
We note that $\sigma$ and $\tau$ are involutions, which is to say that $\sigma\circ\sigma$ and $\tau\circ\tau$ are well defined and equal to the identity on a dense open subset of the plane.  Thus $f$ is birational, and $f^{-1}=\sigma\circ\tau$.  The 2-form $\omega=x^{-1}dx\wedge dy$ is invariant under $f$.
$$
X=\{y=0\},\ \ Y=\{x=0\}, \ \ C=\{(x-1)(y-1)=1\},
$$
we see that $f$ maps $\C^2-(X\cup Y)$ into $\C^2$, and sends $\C^2-(X\cup Y\cup C)$ biholomorphically onto its image.  Our goal in this section will be to construct a complex surface compactifying $\C^2$ in a manner compatible with the dynamics of $f$.

We begin by extending $f$ to a birational map of $\P^2$ via the imbedding $\C^2\ni(x,y)\mapsto[x:y:1]\in\P^2$.  We write the line at infinity as $L_{\infty} =\{[x:y:z]\in\P^{2}:z=0\}$, and thus $\P^{2}=\C^{2}\cup L_{\infty}$.  Using this identification of $\C^{2}$ as a subset of $\P^{2}$, the maps $\sigma$ and $\tau$ may be rewritten in homogeneous coordinates as
$$
\sigma[x:y:z] =[xyz+(-y+z)x^2:xyz+(-x+z)y^2:xyz],
$$
$$
\tau[x:y:z] = [x:b x + (a+1)z-y:z].
$$
It follows that $\tau$ is  invertible and holomorphic on $\P^2$.  On the other hand, $\sigma$ is undefined at the
points of indeterminacy
$$\cI(\sigma)=\cI(f)=\{[0:0:1], [0:1:0], [1:0:0]\}.$$
\medskip
\epsfxsize=1.7in
\centerline{ \epsfbox{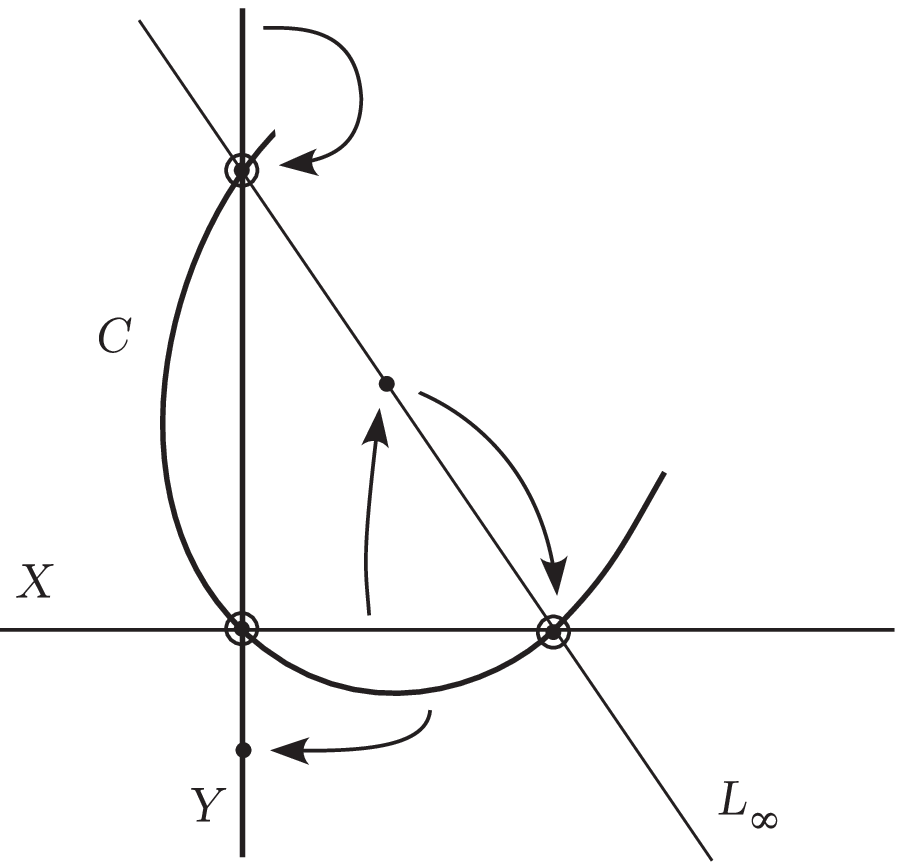}}
\centerline{Figure 1.1.  Behavior of $f$ on $\P^{2}$}
\medskip
If $V$ is a complex curve, then we will denote the {\it proper transform} of $V$ with respect to $f$ and $f^{-1}$ by
$$
f(V) := \overline{f(V-\cI(f)}, \quad f^{-1}(V) := \overline{f^{-1}(V-\cI(f^{-1}))}.
$$
Given a divisor $D$ we let $f^*D$ denote the pullback (i.e. {\it total transform}) of $D$.  This is constructed as follows.  Let $U\subset\P^2-\cI(f)$ and $U'\supset f(U)$ be open sets and $\varphi$ be a meromorphic defining function for $D$ on $U'$.  Then $f^*D$ is defined on $U$ by $\varphi\circ f$.  This specifies $f^*D$ on $\P^2-\cI(f)$.  Extending trivially across $\cI(f)$ completes the construction.  The pushforward of $D$ is simply $f_*D = (f^{-1})^*D$.  We note that if $V$ is an irreducible curve, then $f^*V - f^{-1}(V)$ is always a non-negative linear combination of components of $\cC(f)$.

With these conventions, we have $f(L_\infty) = L_\infty$ and more precisely
$$
f:[x:y:0]\mapsto[x:bx-y:0].
$$
The qualitative behavior of $f$ on $\P^2$ is pictured in Figure 1.1.  The curves $X$, $Y$, and $C$ are all mapped to points and are thus said to be {\it exceptional} curves for the mapping $f$.  The exceptional curves are denoted by heavy lines.  The points of intersection of any two of these exceptional curves are points of indeterminacy.  These points  are marked as dots with circles around them.  The exceptional curves map as follows:
$$
f(C)=(0,a+1)=[0:a+1:1]\in Y\subset\P^2
$$
$$
f(X)=[1:b:0]\in L_\infty, \ \ f[1:b:0]=[1:0:0]=X\cap L_\infty,
$$
$$
f(Y)=[0:1:0]= Y\cap L_\infty
$$
The image points $f(C)$ and $f(X)$ are marked as dots.

We recall the following result from [DF] (also see [FS]).

\proclaim Theorem 1.1. If $f:X\to X$ is a birational self-map of a complex projective surface, then $(f^n)^* D = (f^*)^n D$ for all divisors $D$ and all $n>0$ $\Leftrightarrow$ there is no curve $V\subset X$ such that $f^n(V)\in \cI(f)$ for some $n>0$.  Either of these statements is equivalent to the corresponding statement with $f^{-1}$ replacing $f$.

It was shown in [DF] that the equivalent conditions in Theorem 1.1 can always be arranged to hold by a blowup procedure.
In the present context, Figure 1.1 makes clear that $f(Y)$ and $f^2(X)$ are both points in $\cI(f)$.  So we lift our map $f:\P^2\to\P^2$ to a new manifold $\cX$, which we construct in two steps.  For the first step, we start with $\P^2$ and blow up the three points $B=\{[1:0:0],[1:b:0],[0:1:0]\}\subset L_\infty$.   Blowing up gives us a 2-dimensional manifold $\hat\P^{2}$ and a holomorphic projection $\pi:\hat\P^2\to\P^2$ such that
$\pi:\hat\P^{2}-\pi^{-1}B\to\P^{2}-B$
is a biholomorphic mapping.  We identify curves $\Gamma\subset\P^2$ with their proper transforms $\pi^{-1}(\Gamma) \subset \hat\P^2$.  In Figure 1.2, we let $V_{0}$, $V_{1}$, $C_{0}^{+}$ and $C_{1}^{+}$ denote the proper transforms of $L_{\infty}$, $Y$, $X$, and $C$, respectively.  We denote the exceptional fibers of $\pi$ by dashed lines in Figure 1.2 and label them $V_{2}=\pi^{-1}[0:1:0]$, $V_{5}=\pi^{-1}[1:b:0]$, and $V_{4}=\pi^{-1}[1:0:0]$.  The image of $C_{0}^{+}$ is denoted by a dot in $V_{1}(=Y )$, and the image of $C_{1}^{+}$ is denoted by a dot in $V_{5}$.  We changed the names of the strict transform curves to emphasize the fact that they are now considered in $\hat\P^{2}$ rather than in $\P^{2}$.  For instance, $Y\cap L_{\infty}\ne\emptyset$ in $\P^{2}$, whereas $V_1$ does not intersect $V_0$ in $\hat\P^{2}$.
\medskip
\epsfysize=1.5in
\centerline{ \epsfbox{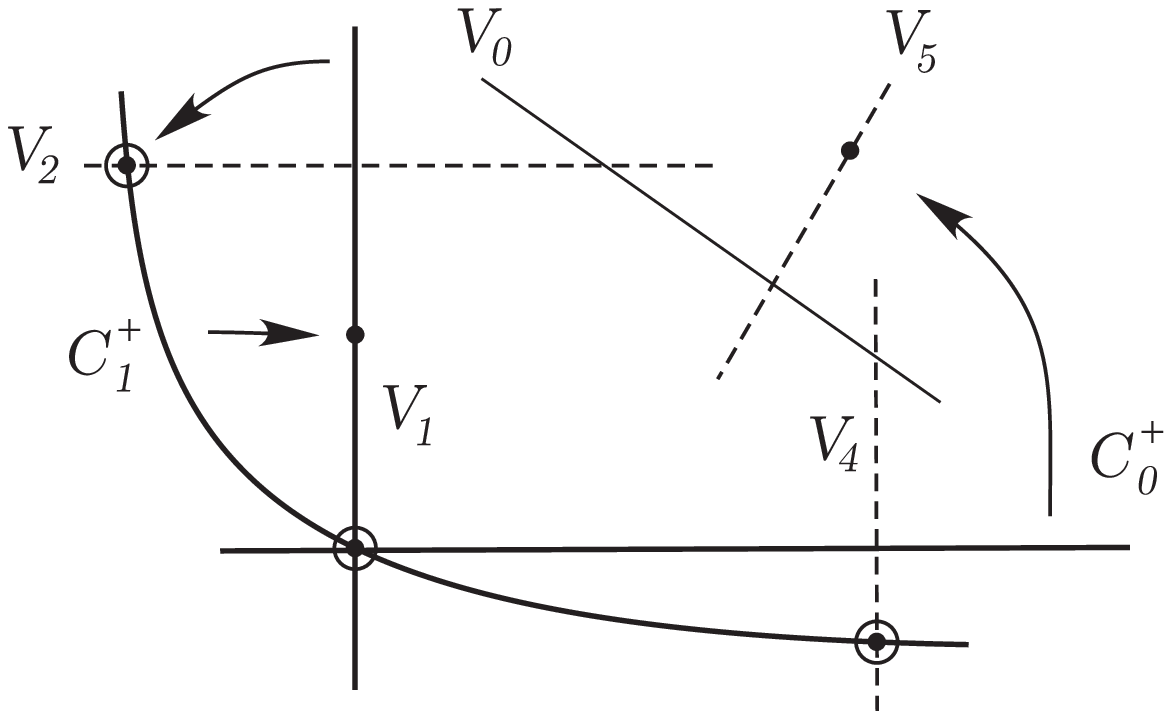}  }
\centerline{Figure 1.2.  Behavior of $f$ on $\hat\P^{2}$}
\medskip

In $\hat\P^{2}$ the curve $V_{1}$ is mapped to the indeterminate point $C^{+}_{1}\cap V_{2}$, so we blow up this point, too. We call the resulting surface $\cX$ and continue to let $\pi:\cX\to\P^2$ denote projection back to $\P^2$.  We let $V_3$ denote the new exceptional fiber of $\pi$.  The manifold $\cX$ is described by Figure 1.3.  The map $f:\cX\to \cX$, has only two exceptional curves $C_{0}^{+}$ and $C_{1}^{+}$ and two points of indeterminacy:
$$
\cI(f) = \{(0,0) = C_0^+\cap C_1^+\cap V_1, \,\,\, C_1^+\cap V_4 = \{y=1\}\cap V_4\}
$$

\epsfysize=1.7in
\centerline{\epsfbox{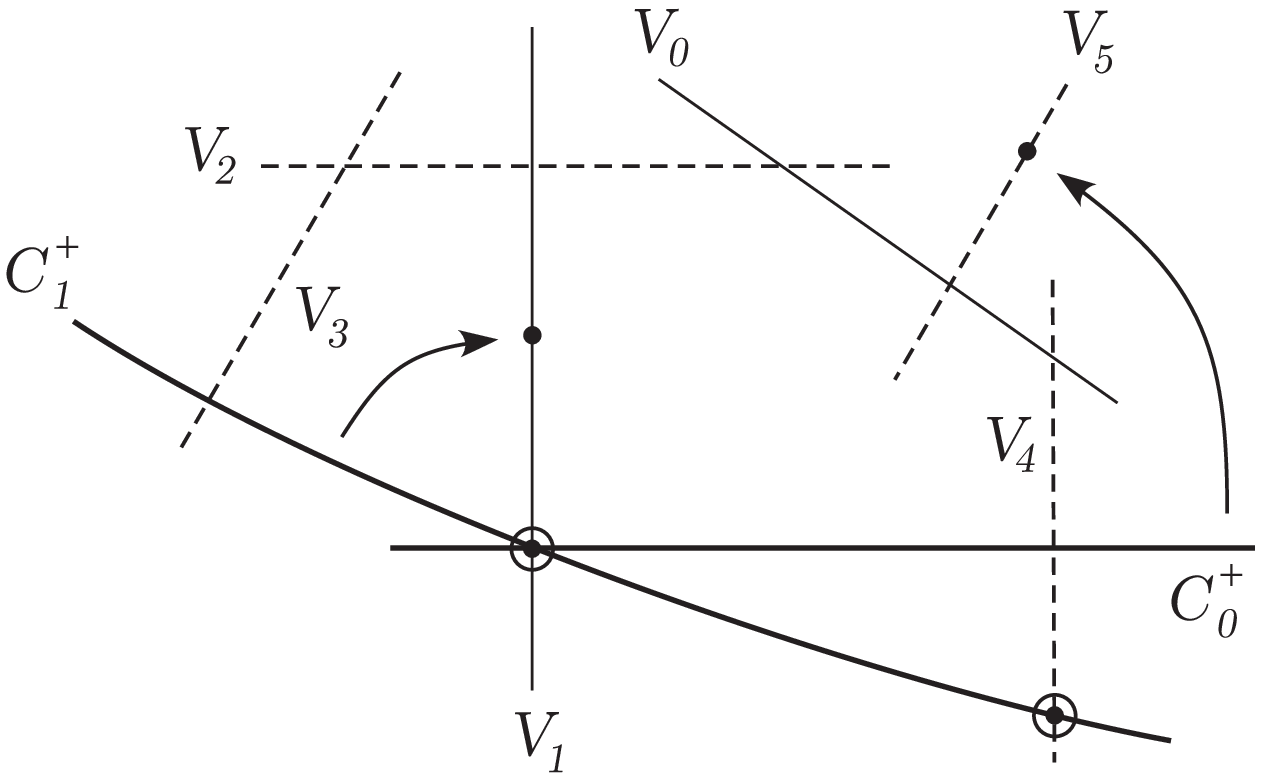}}

\centerline{Figure 1.3.  Behavior of $f$ on $\cX$}
\medskip

To describe the behavior of $f$ on $V_0\cup \dots \cup V_5$, we introduce a complex coordinate $t_{j}$ on each $V_{j}$:

$V_{0}$ is parametrized by its intersection with lines of the form $\{y=t_{0}x\}$.

$V_{1}$ is parametrized by its intersection with lines of the form $\{y=t_{1}\}$.

$V_{2}$ is parametrized by its intersection with lines of the form $\{x = t_{2}\}$.

$V_{3}$ is parametrized by its intersection with conics of the form $\{(y-1)(x-1)=t_{3}\}$.

$V_{4}$ is parametrized by its intersection with lines of the form $\{y=t_{4}\}$.

$V_{5}$ is parametrized by its intersection with lines of the form $\{y=bx +t_{5}\}$.
\smallskip\noindent
Each $V_{j}$, $j\ne1$, consists of points at infinity.  For instance, the point $t_{0}$ in $V_{0}$ is the limit as $s\to\infty$ of the point $(s,st_{0})$.   The point $t_{3}\in V_{3}$ is the limit as $s\to\infty$ of $(1+s^{-1}t_{3},1+s)$, etc.  The action of $f$ on the varieties $V_{j}$ is given as follows:
$$
\eqalign{  V_0\to V_0& {\rm\ \ \ with \ \ }t_0\mapsto t_0=b-t_0,\cr
V_1\to V_3 & {\rm\ \ \ with \ \ } t_1\mapsto t_3=t_1-1,\cr
V_2\to V_2& {\rm\ \ \ with \ \ } t_2\mapsto 1-t_2,\cr
V_3\to V_1& {\rm\ \ \ with \ \ } t_3\mapsto t_1=t_3+a,\cr
V_4\to V_5& {\rm\ \ \ with \ \ } t_4\mapsto t_5=t_4+a,\cr
V_5\to V_4& {\rm\ \ \ with \ \ }t_5\mapsto t_4=t_5+a+1.\cr} \eqno(1.2)
$$

The two exceptional curves are $C_0^+$ and $C_1^+$.   Let us see what happens to $C_0^+$ under $f$.  We find that $f(C_0^+)$ is the point of $V_5$ with $t_5$-coordinate equal to $a$.  Tracking the mapping $V_5\to V_4$, we see that $f^2(C_0^+)$ is the point of $V_4$ with $t_4$-coordinate equal to $2a+1$.  Since the only point of indeterminacy on $V_4\cup V_5$ is the point with $t_4$-coordinate equal to 1, it follows that we can have $f^m(C_0^+)\in\cI$ if and only if $m=2n+2$, and
$$
(n+1)(2a+1) = 1, {\rm\ \ or\ \ } a = -{n\over 2n + 2}.\eqno(1.3)
$$

Next we follow $C_{1}^{+}$.  Since $f(C_1^+)=(0,a+1)\in V_1$, it follows that $f^{2n+1}(C_1^+)\in V_1$ and $f^{2n+2}(C_1^+)\in V_3$.  Since there is no point of indeterminacy in $V_3$, the only way that the orbit of $C_1^+$ can reach $\cI$ is if $f^{2n+1}(C_1^+)=(0,0)$.  Now we know, tracking $V_1\to V_3\to V_1$ in the $t_1$ and $t_3$ coordinates above, that  $f^2$ maps $V_1$ to itself by $(0,y)\mapsto (0,y+a-1)$.  Thus we have $f^{2n+1}(C_1^+)\in\cI$ exactly when
$$a+1+n(a-1)=0, {\rm \ \ or \ \ }a={n-1\over n+1}.\eqno(1.4)$$

If $a$ satisfies neither (1.3) nor (1.4), then the forward orbit of an exceptional curve never encounters a point of indeterminacy. So by Theorem 1.1, $(f^n)^* = (f^*)^n$ for all $n\in\Z$.

\bigskip\noindent{\bf \S 2.  Action of the Pullback on Divisors. }  We say that two divisors $D_{1}$ and $D_{2}$ are {\it linearly equivalent} if $D_{1}-D_{2}$ is the divisor of a rational function.  We let ${\rm Pic}({\cal X})$ (the Picard group) denote the set of divisors on ${\cal X}$ modulo linear equivalence.  In ${\bf P}^{2}$ all lines are equivalent, and ${\rm Pic}({\bf P}^{2})$ is generated by the equivalence class of a line.  The operation of blowing up a point increases the dimension of ${\rm Pic}(\cX)$ by one, so ${\rm Pic}({\cal X})$ has dimension 5.   In fact, ${\cal V}:=\{V_{j}: 1\le j\le 5\}$ defines a basis of ${\rm Pic}({\cal X})$.

The action of $f$ on the various curves $V_{j}$ is summarized in (1.2) and Figure 1.3. Since both $C_0^+$ and $C_1^+$ are critical with multiplicity one, we see:
$$
\eqalign{ f^{*}: \ \ & V_{0}\mapsto V_{0},\ \ \ \  V_{4}\mapsto V_{5}\mapsto V_{4}+C_{0}^{+}\cr
&V_{2}\mapsto V_{2},\  \ \ \  V_{3}\mapsto V_{1}\mapsto V_{3} + C_{1}^{+}.\cr}\eqno(2.1)
$$

Pullback preserves linear equivalence among divisors, so (2.1) serves to define $f^{*}$ as a linear map of the Picard group. In order to find the matrix for $f^*$ relative to ${\cal V}$, let us consider the pullback map $\pi^{*}:{\rm Pic}({\bf P}^{2})\to {\rm Pic}({\cal X})$.  We let $L\in{\rm Pic}(\cX)$ be the class obtained by pulling back a generic line in $\P^2$.  In order to  determine $\pi^{*}C_{1}^{+}$, we consider the composition $\pi=\pi_{1}\circ\pi_{2}$, with $\pi_{1}:\hat{\bf P}^{2}\to{\bf P}^{2}$ and $\pi_{2}:{\cal X}\to\hat{\bf P}^{2}$.  For any curve $C$, the difference $\pi^*C -\pi^{-1}(C)$ is a sum of exceptional fibers that map to $C$.  Hence
$$
\pi^{*}_{1}C^{+}_{1}=C^{+}_{1}+V_{2}+V_{4}\in{\rm Pic}(\hat{\bf P}^{2}).
$$
Likewise
$$
\pi^{*}C^{+}_{1}= \pi^{*}_{2}(\pi^{*}_{1}C^{+}_{1}) = \pi^{*}_{2}C_{1}^{+}+\pi^{*}_{2}V_{2}+\pi^{*}_{2}V_{4}
= (C^{+}_{1}+V_{3}) + (V_{2}+V_{3}) +V_{4}.
$$
Because $C^+_1$ and $2L$ are linearly equivalent in $\P^2$, we conclude that
$$
2L = C^{+}_{1}+V_{2}+2V_{3}+V_{4} \eqno(2.2)
$$
Similarly, we find
$$
L=\pi^{*}V_{0}=V_{0}+V_{2}+V_{3}+V_{4}+V_{5},
$$
$$
L=\pi^{*}V_{1}=V_{1}+V_{2}+V_{3},\eqno(2.3)
$$
$$
L=\pi^{*}C_{0}^{+}=C_{0}^{+}+V_{4}.
$$

Equation (2.2) and the last two equations in (2.3) allow us to solve for $C_0^+$ and $C_1^+$ in terms of the basis
${\cal V}$.  Then from (2.1--2.3) we can find the matrix for $f^*$ relative to ${\cal V}$.
The result is
$$
C_0^+ = \pmatrix{ 1\cr 1\cr 1\cr -1\cr 0\cr}, \quad C_1^+ =\pmatrix{2\cr 1\cr 0\cr -1\cr 0}, \quad
f^* = \pmatrix{ 2&0&1&0&1\cr
1&1&0&0&1\cr
1&0&0&0&1\cr
-1&0&0&0&0\cr
0&0&0&1&0\cr}.\eqno(2.4)$$
The characteristic polynomial of $f^*$ is
$$
(x-1)^{2}(x^{3}-x^{2}-2x-1).\eqno(2.5)
$$
\proclaim Theorem 2.1.  Let $a,b\in{\bf C}$, $b\ne0$ be given, and suppose that (1.3) and (1.4) do not hold.  Then the asymptotic degree growth rate for $f$ is given by the largest root $\rho\sim 2.1479$ of the polynomial (2.5).

\noindent{\it Proof.}  Let $\|\cdot\|$ denote any norm on ${\rm Pic}(\cX)$ and $\|f^{n*}\|$ be the corresponding operator norm of $f^{n*}$.  By Theorem 1.1 and our hypothesis on $a$ and $b$, we have $(f^n)^*=(f^*)^n$.  Hence the
quantity
$$
\lim_{n\to\infty}\|f^{n*}\|^{1/n}
$$
is equal to the spectral radius $\rho$ of the matrix (2.4). On the other hand, this quantity is independent of birational changes of coordinate, and on $\P^2$ it is equal to the growth rate of the degree of $f^n$.
\medskip

Now we consider the intersection product on ${\rm Pic}({\cal X})$.
The intersection number (see [GH, pp.\ 470-476]) of distinct irreducible curves, written $D_{1}\cdot D_{2}$, is simply the number of set-theoretic intersections counted with multiplicity.  The self-intersection $D\cdot D$ of a line in ${\bf P}^{2}$ is $+1$, and the self-intersection of the exceptional fiber of a blowup is $-1$.  The self-intersection of a curve $D$ decreases by one every time we blow up a smooth point on $D$.  The intersection number depends only on the classes in ${\rm Pic}({\cal X})$ of the divisors concerned.  Relative to ${\cal V}$, the intersection form for ${\rm Pic}({\cal X})$ has matrix
$$\pmatrix{ 0&1&0&0&0\cr
1&-2&1&0&0\cr
0&1&-1&0&0\cr
0&0&0&-1&0\cr
0&0&0&0&-1\cr}.\eqno(2.6)$$

The curves $V_0$ and $V_1$ are invariant under $f^*$, and one can check that the same is true of the (intersection) orthogonal complement $S=\{V_0,V_1\}^\perp$.  Hence $f^*|S$ has irreducible characteristic polynomial $x^{3}-x^{2}-2x-1$.
As it happens, $S$ is spanned by $C_0^+$, $C_1^+$ and $f^* C_0^+$, and it will be particularly convenient to perform subsequent computations using these three curves as a basis.  Relative to this basis, $f^{*}$ has the form
$$
\pmatrix{ 0&1&1\cr
0&0&1\cr
1&1&1\cr}. \eqno(2.7)
$$
Recall that $f$ and $f^{-1}$ are conjugate via $\tau$ and that in particular, $C_j^- = \tau(C_j^+)$ and $\tau(V_j) = V_j$ for $j=0,1$.  Hence $S$ is also spanned by $C_0^-$, $C_1^-$, $f_*C_0^-$, and with respect to this alternative basis $f_*$, too, has matrix (2.7).  Finally, let us give the intersection matrix relative to these two bases.  If $W=a_{1}C_{0}^{+}+a_{2}C_{1}^{+}+a_{3}f^{*}C_{0}^{+}$ and $W'=b_{1}C_{0}^{-}+b_{2}C_{1}^{-}+b_{3}f_{*}C_{0}^{-}$, then
$$
W\cdot W' =(a_{1},a_{2},a_{3})\pmatrix{ 1&2&2\cr
2&2&3\cr
2&3&5\cr} \pmatrix{b_{1}\cr b_{2}\cr b_{3}\cr}.\eqno(2.8)
$$

\bigskip\noindent{\bf \S 3.  Real Mappings and their Compactifications. }
Here we review the results of the construction of ${\cal X}$ in \S1, now restricting to real parameters, and considering the action of $f$ on real points.    We consider $\R^2$ as a subset of $\P^2$, using the same imbedding $\R^2\ni(x,y)\mapsto[x:y:1]\in\P^2$ as before.  Let $\overline{\R^2}$ denote the closure of $\R^2$ in $\P^2$.  It follows that the intersection with the line at infinity is the circle $\overline{\R^2}\cap L_\infty=\{[x:y:0]=[-x:-y:0]: x,y\in\R\}$, and $\overline{\R^2}$ is a real analytic submanifold of ${\bf P}^{2}$ which is naturally identified with the real projective plane.  We have pictured the real projective plane on the left hand side of Figure 3.2; it is represented as a closed disk with opposite sides identified in the fashion indicated by the arrows.  The boundary of the disk is the circle $V_{0}\cap {\bf R}^{2}$, which we have drawn as a hexagon so that each side contains one of the points that must be blown up to obtain $\cX$.  In addition, we have included the intersection of $\overline{\R^2}$ with the sets $X$, $Y$, $C$ from Figure 1.1, as well as three points at infinity.
\medskip
\epsfysize=.9in
\centerline{ \epsfbox{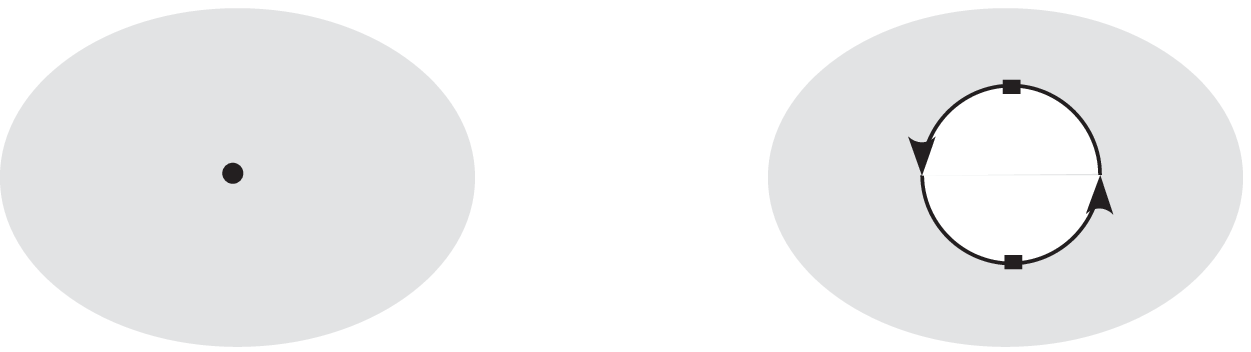}  }
\centerline{Figure 3.1.  A real (unoriented) blow-up.}
\medskip

Now let ${\cal X}_{R}$ denote the closure of ${\bf R}^{2}$ inside ${\cal X}$.  The complex manifold  ${\cal X}$ was constructed by blowing up four points over  $\overline{\R^2}\cap L_\infty$.  We recall that blowing up a point in a real surface amounts topologically to replacing the point with a cross-cap, as pictured in Figure 3.1.  Take $[1:b:0]$, for instance.  Blowing up this point has the effect of `inserting' two segments into the hexagon on the left side of Figure 3.2.  These appear in the right side of the figure labeled as $V_5$.  Note that the identification between them is opposite that of the surrounding points in $V_0\cap\cX_R$.  Each of the other blowups likewise modifies the hexagon: opposing pairs of segments labeled $V_j$, $j=2,3,4$ are also added at infinity and identified in a manner counter to that of the surrounding points.  Since adding two cross-caps is the same as adding a handle, we find that ${\cal X}_{R}$ is homeomorphic to the real projective plane with two handles attached, and the first homology group is $H_{1}({\cal X}_{R};{\bf Z})\cong{\bf Z}^{4}\oplus({\bf Z}/2{\bf Z})$.

\medskip
\epsfysize=2in
\centerline{ \epsfbox{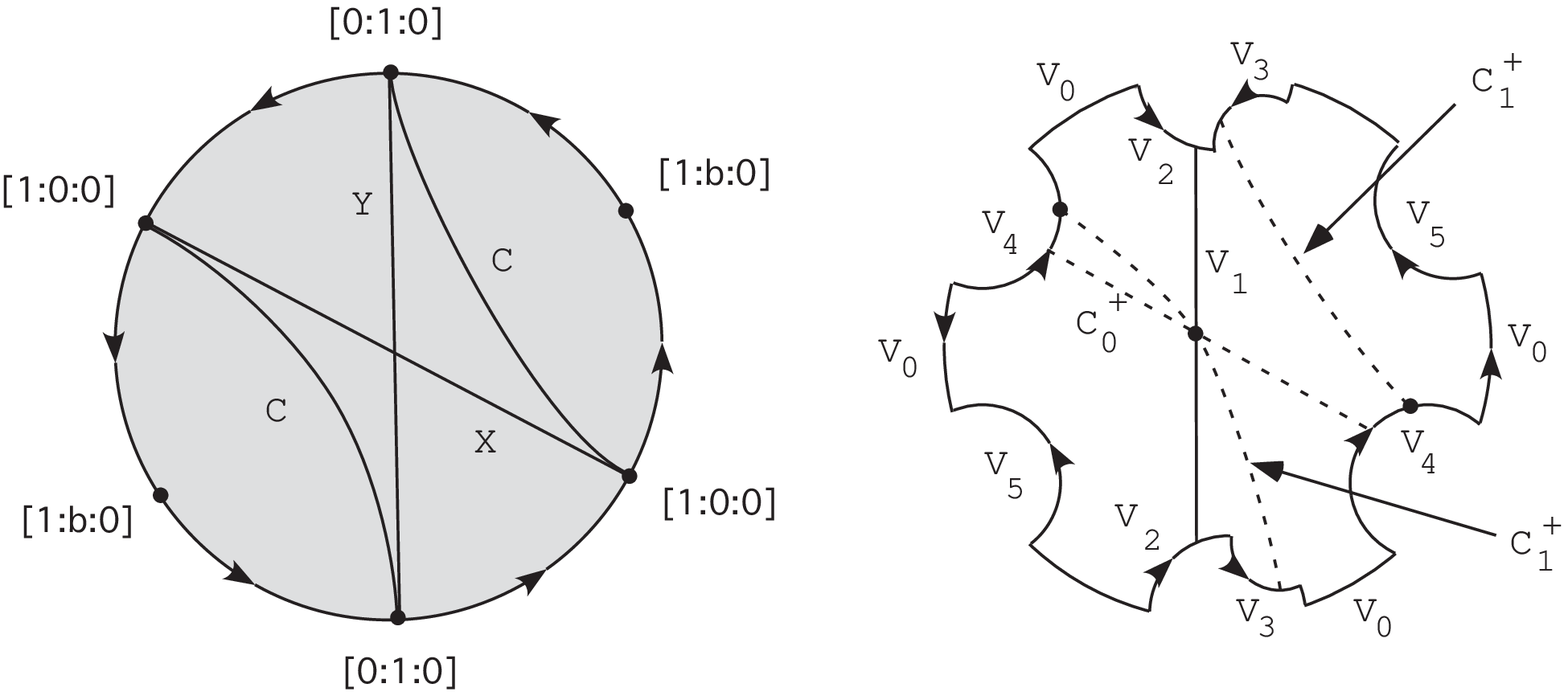}  }
\centerline{Figure 3.2.  Manifold ${\cal X}_{R}$.}
\medskip
If $a,b$ are real, then $f$ defines a birational map $f_{R}$ of the surface ${\cal X}_{R}$.  The points of indeterminacy of $f$ and $f_{R}$ are essentially the same.  They  are the origin $(0,0)\in V_{1}$ and $C_{1}^{+}\cap V_{4}=\{t_{4}=1\}$.  The involution $\tau$ flips points about the horizontal axes of the hexagons in Figure 3.2.

The sets ${V_0\cap \bf R}^2$ and ${V_2\cap \bf R}^2$ will play an unimportant role in our subsequent discussion.  We remove them and rewrite $\cX_{R}-(V_{0}\cup V_{2})$ as the ``dynamic hexagon'' in Figure 3.3.  The vertices of this hexagon represent the points of $V_0\cup V_2$ which were removed.  The opposite sides of the hexagon are identified as in Figure 3.2.  The dotted arrows indicate the direction in which points in $V_{j}$, $j=1,3,4,5$, move under $f^2$ for the parameter range $a<-1,b>0$.  The critical locus and its forward image are also shown. In the sequel we will see how dynamical behavior may be deduced from this drawing.
\medskip
\epsfysize=1.9in
\centerline{ \epsfbox{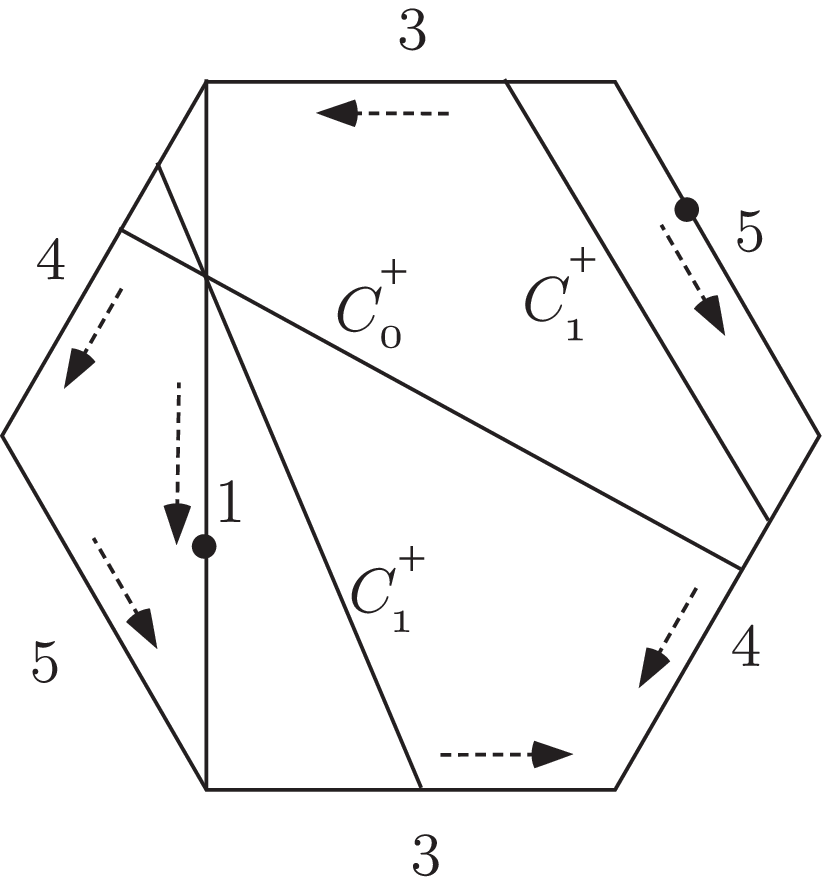}  }
\centerline{Figure 3.3.  Fundamental hexagon.}
\medskip

The construction of $\cX_{R}$ in the case $b=0$ is slightly different.  In this case, the picture on the left half of Figure 3.2 will be a rectangle, since the sides corresponding to $[1:b:0]$ and $[1:0:0]$ will be the same.  The picture on the right hand half will similarly be altered by the elimination of $V_{5}$.  In [BHM2] it was shown that $f_{a,b}$ is birationally conjugate to the map $g_{a}=(y(x+a)/(x-1),x+a-1)$ when $b=0$.  This simpler form was used to study $f_{a,0}$ in [A1--7] and in [BD1].  The torus compactification of ${\bf R}^{2}$ is natural for $g_{a}$ for generic $a$.  It is interesting to compare the differences between the stable/unstable laminations of $f_{-2,0}$ in ${\cal X}_R$, as shown in Figure 0.2,  and in the torus compactification, as shown in Figure 9 of [BD1].

\bigskip\noindent{\bf \S4.  Combinatorics of Real Curves. }
In this section we begin our discussion of the dynamics of $f$ for real parameters $a$, $b$.  The goal will be to find families of real curves on which the action of the real map $f^{-1}:\cX_R\to\cX_R$ mimics that of the pullback action
$f^*:{\rm Pic}(\cX)\to{\rm Pic}(\cX)$.  Our analysis begins with the invariant curves $V_j$, $j=1,3,4,5$.  While these curves are ``properly invariant'' in the sense that $f^2(V_j) = V_j$, they are not ``totally invariant'' in the sense of pushforward.  Rather, we have
$$
f_*:V_3\mapsto V_1\mapsto V_3 + {\cal C}^-_1, \qquad V_5\mapsto V_4\mapsto V_5+{\cal C}^-_0.\eqno(4.2)
$$
In contrast, the curves $V_0$ and $V_2$ are invariant in both senses and will play no role in our analysis.  For each $j=1,3,4,5$, we declare $E_j^s,E_j^u\subset V_j\cap{\cal X}_R$ to be the smallest closed intervals such that
$$
E^s_j\supset V_j\cap\bigcup_{n\ge0}f^n{\cal C}(f^{-1}),\qquad  E^u_j\supset V_j\cap\bigcup_{n\ge0}f^{-n}{\cal C}(f).  \eqno(4.3)
$$
In particular, $f^{-1}(\bigcup E_j^u) \subset \bigcup E_j^u$ and $f^{-1}(\cC(f^{-1})) = \cI(f)\subset \bigcup E_j^u$.  Hence $f$ restricts to an everywhere well-defined differentiable map of the open set $\cX_R-\bigcup E^u_j$ into the smaller region $\cX_R-\bigcup E^u_j-\cC(f^{-1})$.  The inverse of $f$ is similarly well-behaved on $\cX_R-\bigcup E^s_j$.
\medskip
\epsfysize=2in
\centerline{ \epsfbox{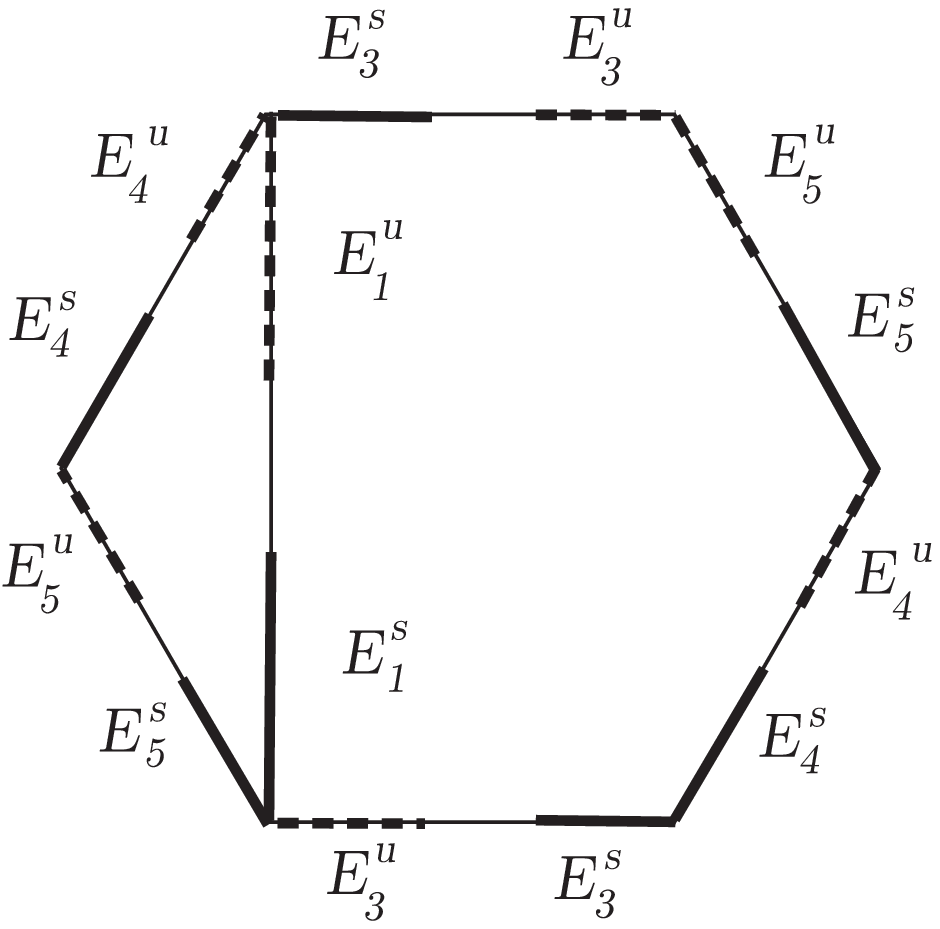}  }
\centerline{Figure 4.1.  The sets $E^{s/u}_{j}$.}
\bigskip

For the rest of this paper we will assume that
$$
a<-1,\ \  b\ne0.\eqno(4.4)
$$
Under this assumption, and in terms of the parametrizations for $V_1,V_3,V_4,V_5$ given in \S2, direct computation
reveals that
$$
\eqalign{E_1^s &= [-\infty,a+1],\cr E_1^u &= [0,\infty],}\quad
\eqalign{E_3^s &= [-\infty,-1], \cr E_3^u &= [1,\infty],}\quad
\eqalign{E_4^s &= [-\infty,2a+2], \cr E_4^u &= [0,\infty],}\quad
\eqalign{E_5^s &= [-\infty,a+1], \cr E_5^u &= [-1-a,\infty].}
$$
\centerline{\epsfysize=2in \epsfbox{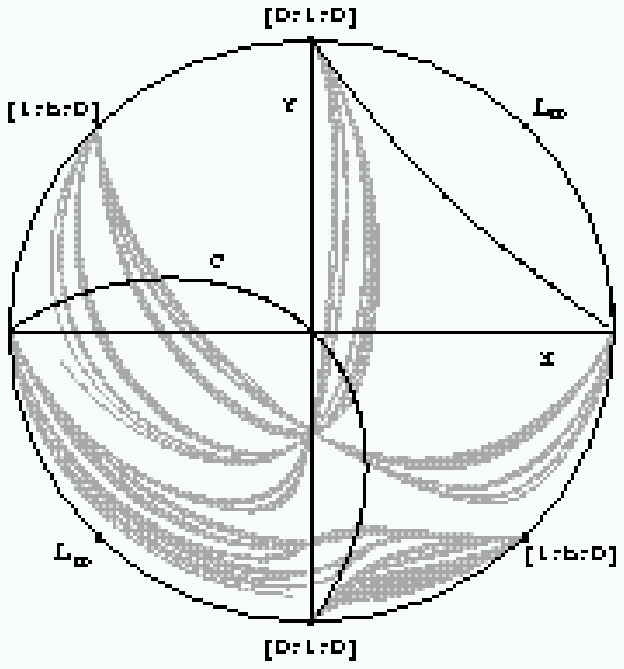} \hfil\epsfysize=2in \epsfbox{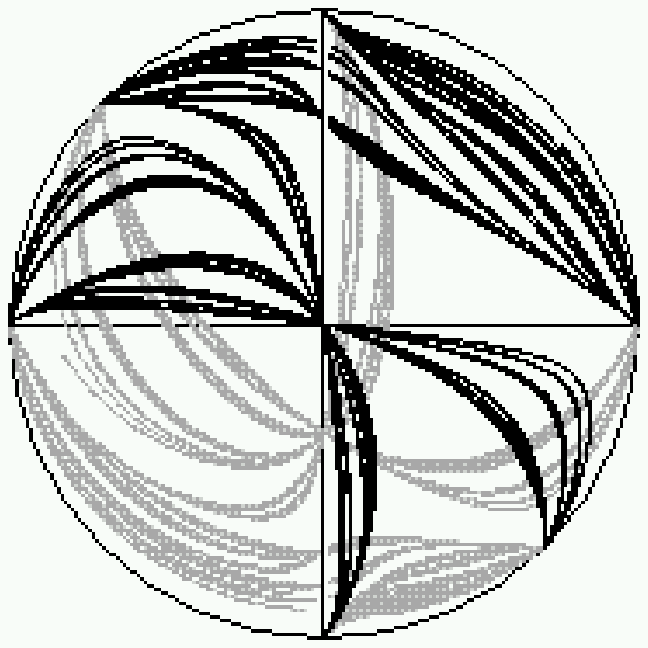}  }
\centerline{Figure 4.2.  $f_{a,b}$ with $a=-2$, $b=-1$; an unstable manifold (left),}

\centerline{stable and unstable manifolds (right).}
\medskip

Figure 4.1 illustrates the arrangement of the intervals, superimposing them on the hexagon from Figure 3.3.
It is important for the analysis to follow that $E^s_j$ and $E^u_j$ be disjoint for each $j$.  The reason for the restriction
(4.4) is that when $a>-1$, there always exists an index $j$ for which the overlap $E^s_j\cap E^u_j$ is non-trivial.  All of our arguments in the sequel will be given in the case $b>0$; we leave to the reader the modifications for the case $b<0$.  In Figure 4.2 we present an unstable manifold with $b<0$, together with notation from Figure 1.1.  All other figures in this paper are drawn for $b>0$.  For the case $b<0$ Figure 3.3, for example, would be re-drawn with the locations of $V_4$ and $V_5$ interchanged.
To describe the action of $f$, we let $R^\pm_j$, $j=1,\dots,7$, denote the connected components of
$\cX_R -\cC(f^{\pm 1}) - \bigcup V_j$ (see Figure 4.3). Since $f$ maps $\cX_R - \cC(f) - \bigcup V_j$ diffeomorphically onto
$\cX_R-\cC(f^{-1}) - \bigcup V_j$, we see that for each $j$ there is a $k$ such that $fR^+_j=R^-_k$. Specifically, we have
$$
\eqalign{ f(R^+_1)=R^-_2,&\qquad f(R^+_5)=R^-_7\cr
f(R^+_2)=R^-_1,&\qquad  f(R^+_6)=R^-_6\cr
f(R^+_3)=R^-_3,&\qquad  f(R^+_7)=R^-_5\cr
f(R^+_4)=R^-_4.&\cr}\eqno(4.5)
$$

To each region $R_j^+$ we associate the family of {\it $j$s-arcs}, which join suitable pairs of intervals $E_j^u$ on
the boundary of $R_j^+$.  For example, we declare an arc $\gamma$ to be a 1s-arc if it lies to the left of $V_0$, and if
one of its endpoints is in $E^{u}_{1}\cup E^{u}_{4}$ and the other endpoint is in $E^{u}_{5}$. We define  $j$s-arcs for $j=2,3,4$ by referring to Figure 4.4. The general idea is that each family of $js$-arcs should contain an arc $\gamma$ in $\cC(f)$ and further include all arcs homotopic to $\gamma$ through deformations that leave endpoints in $\bigcup E_j^u$ and interior points in $\cX_R-\bigcup V_j$.  We similarly define (referring again to Figure 4.4) $j$u-arcs for $j=1,2,3,4$ joining pairs of intervals $E_j^s$.

\epsfysize=1.8in
\centerline{ \epsfbox{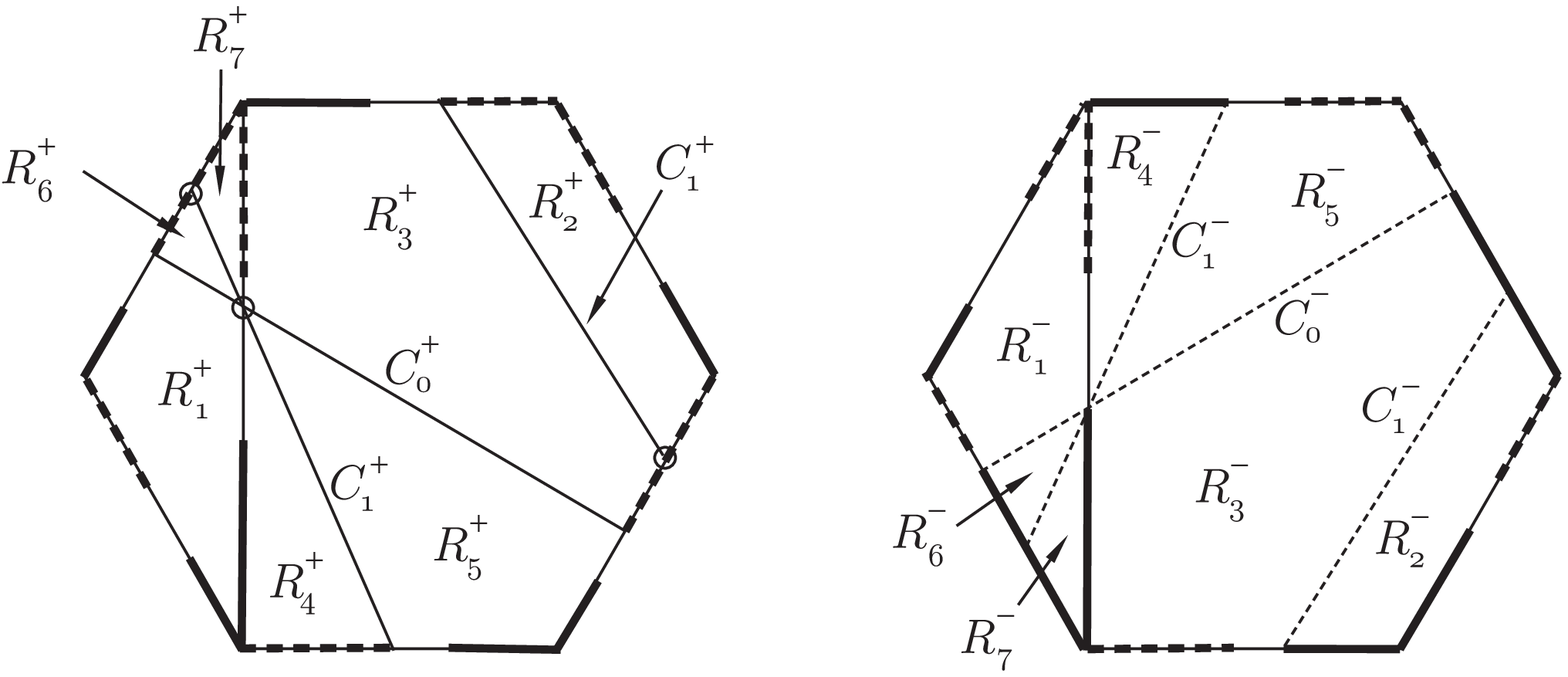}  }
\centerline{Figure 4.3.  Partition of $\cX_\R$ into sets $R^{+}_{j}$ (on left) and $R^{-}_{j}$ (on right).}

\medskip
\epsfysize=1.5in
\centerline{ \epsfbox{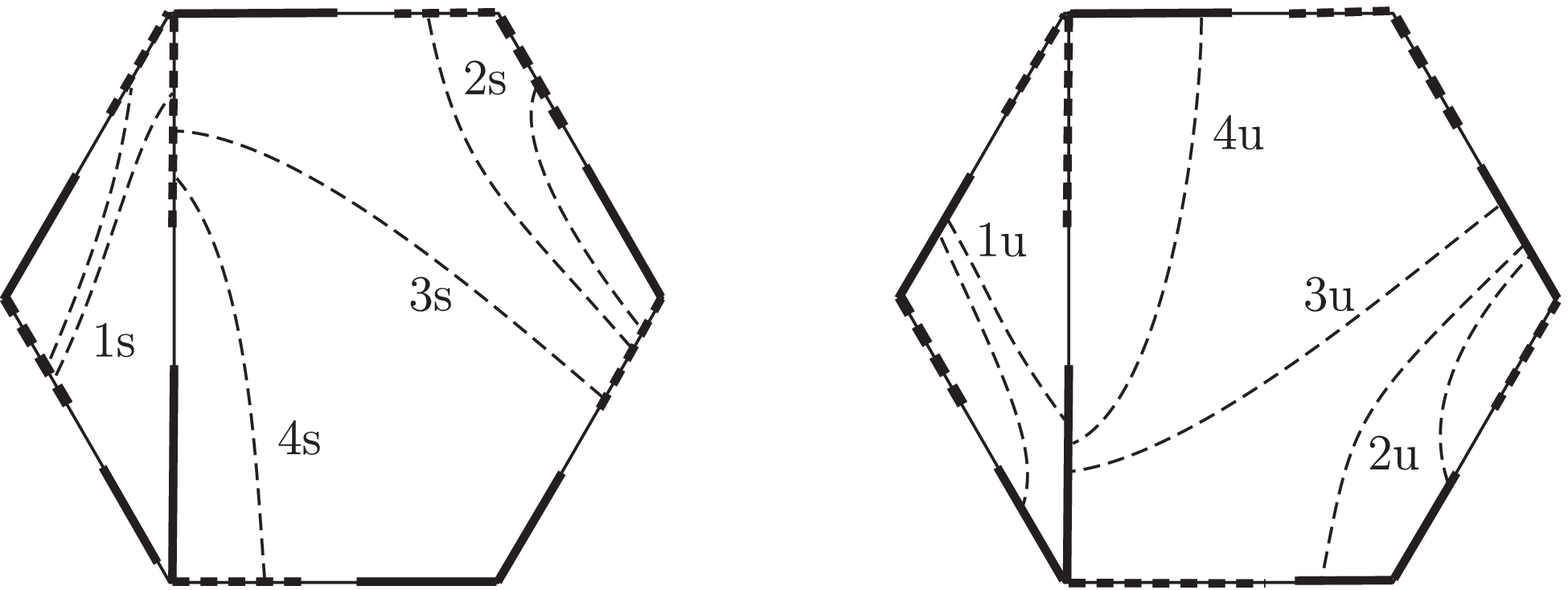}  }
\centerline{Figure 4.4.  s-arcs and u-arcs.}
\medskip

With the aid of Figure 4.3, it is a routine matter to determine the preimage of a $j$s-arc.  For instance, the
left half of Figure 4.5 shows a 4s-arc $\gamma$ drawn relative to the regions $R^{-}_{k}$.  Thus we see $\gamma$ must cross $R^{-}_{4}$, $R^{-}_{5}$ and $R^{-}_{3}$.  From (4.5) and the facts that $f^{-1} C^{-}_{1}=(0,0)$ and $f^{-1} C^{-}_{0}=\{t_{4}=1\}\in V_{4}$, we conclude that $f^{-1}\gamma$ must cross $R^{+}_{4}$, $R^{+}_{7}$ and $R^{+}_{3}$ as shown on the right side of Figure 4.5.  In particular, $f^{-1}\gamma$ contains sub-arcs of type 3s and 4s.  If we repeat this reasoning for all the arcs of type $j$s, we find:

\proclaim Proposition 4.1.  The s-arcs are mapped under $f^{-1}$ as follows:
\vskip0pt  If $\gamma$ is a 1s-arc, then $f^{-1}\gamma$ contains a 2s-arc;
\vskip0pt  If $\gamma$ is a 2s-arc, then $f^{-1}\gamma$ contains arcs of type 1s and 3s;
\vskip0pt  If $\gamma$ is a 3s-arc, then $f^{-1}\gamma$ contains arcs of type 1s, 3s, and 4s;
\vskip0pt  If $\gamma$ is a 4s-arc, then $f^{-1}\gamma$ contains arcs of type 3s and 4s.
\vskip0pt\noindent Further, each of the $j$s-subarcs of $f^{-1}\gamma$ is contained in $R^{+}_{j}$.

Similar reasoning shows that Proposition 4.1 remains valid if we replace ``$f^{-1}$'' by ``$f$'' and ``$j$s-arc'' by ``$j$u-arc.''

If $\gamma$ is an arc in $\cX_{R}$, we use the notation ${\bf s}(\gamma)=(n_{1},n_{2},n_{3},n_{4})$ to mean that $\gamma$ contains $n_{j}$ distinct subarcs of type $j$s.  We interpret the word `distinct' very strongly here, saying that two s-arcs of a given type are distinct if they meet (if at all) only at endpoints. Similarly, we use the notation ${\bf u}(\gamma)=(m_{1},m_{2},m_{3},m_{4})$ to mean that $\gamma$ contains $m_{j}$ distinct sub-arcs of type $j$u.

In light of Proposition 4.1 we define the matrix
$$F=\pmatrix{0&1&1&0\cr
1&0&0&0\cr
0&1&1&1\cr
0&0&1&1\cr}.
\eqno(4.6)
$$
Note that if $\gamma_1$ and $\gamma_2$ are distinct s-arcs, then the s-arcs in $f^{-1}(\gamma_1)$ and $f^{-1}(\gamma_2)$
will all be distinct because $f^{-1}$ is a diffeomorphism off $\cC(f^{-1})$, and $f^{-1}(\cC(f^{-1}))$ lies on the boundaries
of the various regions $R_j^-$.  Therefore, we immediately arrive at

\proclaim Corollary 4.2.  If $\gamma$ is an arc in $\cX_{R}$, then
${\bf s}(f^{-1}\gamma)^{t}\ge F{\bf s}(\gamma)^{t}$ and ${\bf u}(f^{-1}\gamma)^{t}\ge F{\bf u}(\gamma)^{t}$.

\epsfysize=1.9in
\centerline{ \epsfbox{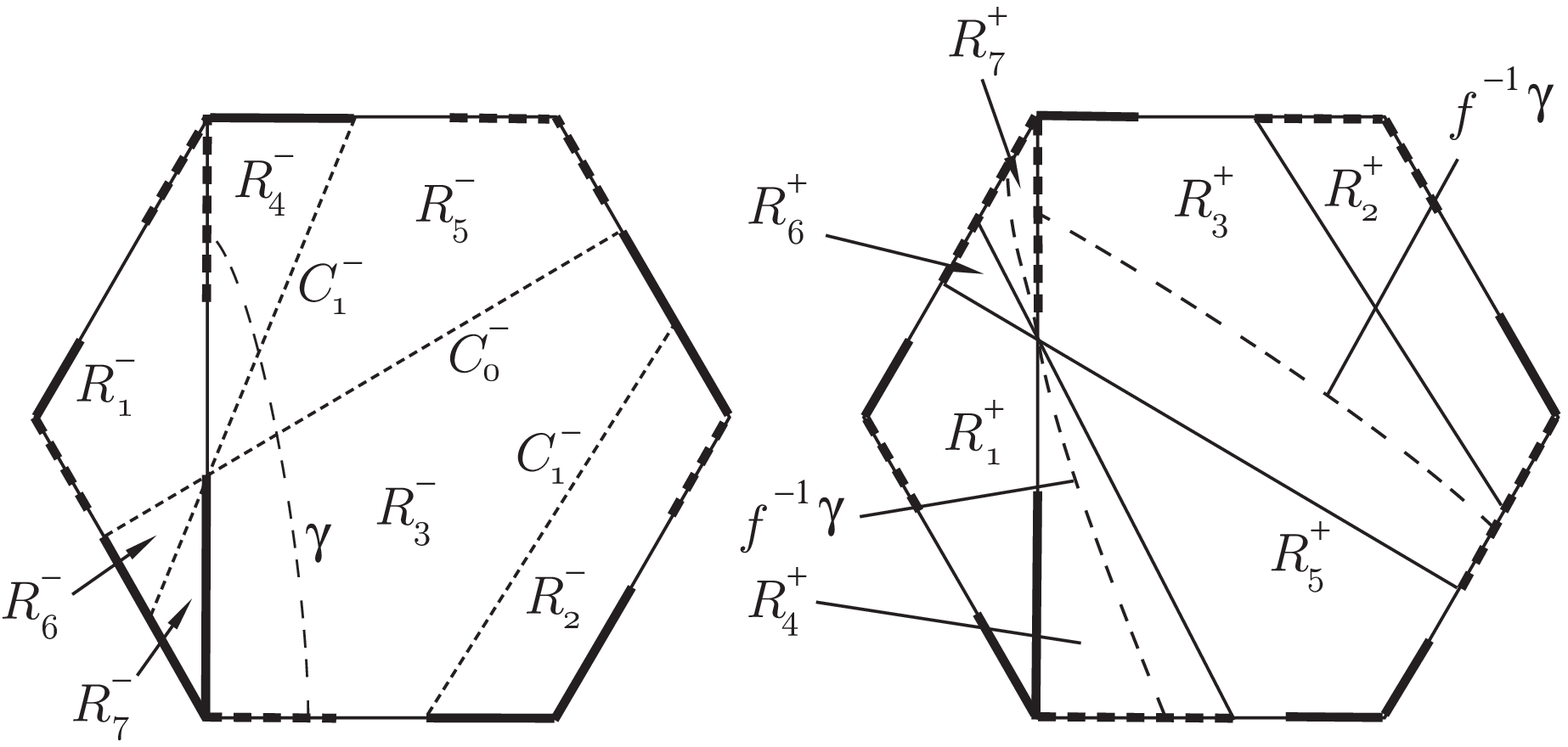}  }
\centerline{\rm Figure 4.5.  Preimage of a 4s-arc.}
\bigskip

We define
$$
Q =
\pmatrix{1&0&0&0\cr
0&1&0&1\cr
0&0&1&1\cr
0&1&1&1\cr},
\eqno(4.7)
$$
signifying that if $(i,j)$ is a pair for which $Q_{i,j}=1$, then every $i$s-arc must meet every $j$u-arc in at least one point.  Since an endpoint of an s-arc cannot meet a u-arc at an endpoint, we immediately conclude

\proclaim Proposition 4.3.  If $\gamma_{1}$ and $\gamma_{2}$ are arcs in $\cX_{R}$, then $\gamma_{1}$ must intersect $\gamma_{2}$ in at least ${\bf s}(\gamma_{1})Q{\bf u}(\gamma_{2})^{t}$ points.

\proclaim Theorem 4.4.  Let $a<-1$, $b\ne0$, and let $\gamma$ be a curve that is disjoint from $\bigcup_{n\ge0}f^{-n}\cI(f)$ and that contains a u-arc.  Then ${\bf u}(f^n\gamma)$ and ${\rm Length}(f^n\gamma)$ both grow at least as fast as $\rho^n$  as $n\to\infty$.  If $\gamma'$ is a curve that is disjoint from $\bigcup_{n\ge0}f^n\cI(f^{-1})$, then the number of intersection points  $\#(f^n\gamma\cap f^{-n}\gamma')$ grow at least as fast as $\rho^{2n}$ as $n\to\infty$.

\noindent{\it Proof. }  By Corollary 4.2, the number of u-arcs in $f^{n}\gamma$ is at least $F^{n}{\bf u}(\gamma)^{t}$.  The characteristic polynomial of $F$ is $(x-1)(x^{3}-x^{2}-2x-1)$, so its largest eigenvalue is $\rho$.  Moreover, all entries of the matrix $F^3$ are positive.  Therefore by the Perron-Frobenius theorem, the vector ${\bf u}(\gamma)$, whose entries
are negative, will be attracted to the eigenspace of $\rho$ under iteration of $F$.  In particular the growth of ${\bf u}(f^{n}\gamma)$ is comparable to $\rho^{n}$.  Since the minimal length of a $u$-arc is bounded away from zero, it follows
that the length of $f^n\gamma$ also grows at least as quickly as $\rho^n$.  Similarly, ${\bf s}(f^{-n}\gamma')$ grows like $\rho^{n}$, and the second statement follows from Proposition 4.3.
\bigskip

The connection between the previous several results and those of \S2 can be seen
more clearly by expressing everything relative to the curves $C_0^+, C_1^+, f^*C_0^+$.

\proclaim Proposition 4.5.   We have
$$
\eqalign{
{\bf s}( C^+_0)={\bf u}( C^-_0)& =(0,0,1,0),  \cr
{\bf s}( C^+_1)={\bf u}( C^-_1)& =(0,1,0,1),  \cr
{\bf s}(f^* C^+_0)={\bf u}(f_* C^-_0)& =(1,0,1,1).\cr}
$$

\noindent{\it Proof.}  The first two lines are seen by referring to Figures 4.3 and 4.4.  By Corollary 4.2, both ${\bf s}(f^{*} C_{0}^{+})$ and ${\bf u}(f_{*} C^{-}_{0})$ dominate the vector $(1,0,1,1)$.  That is, for example,
$$
{\bf s}(f^{*} C_{0}^{+})=(1,0,1,1)+(v_{1},v_{2},v_{3},v_{4})
$$
for a 4-tuple of integers $v_{j}\ge0$.  Since the entries of $Q$ are nonnegative, we have
$$
{\bf s}(f^{*} C_{0}^{+})\, Q\,{\bf u}(f_{*} C^{-}_{0})^{t}\ge((1,0,1,1)+{\bf v})\,Q(1,0,1,1)^{t}=5+v_{1}+v_{2}+2v_{3}+2v_{4}.
$$
By Proposition 4.3, this number is a lower bound for the number of real intersections between the curves $f^{-1}( C^{+}_{0}\cap{\bf R}^{2})$ and $f( C^{-}_{0}\cap{\bf R}^{2})$.

On the other hand, we may consider the complex intersection number
$$
f^{*} C^{+}_{0}\cdot f_{*} C^{-}_{0}= f^{*2} C^+_0\cdot C^-_0=5
$$
by (2.4) and (2.6).  Since the number of real intersections is no larger than the number of complex intersections we have ${\bf v}=0$, which completes the proof.
\bigskip

From Proposition 4.5, one may compute directly that

\proclaim Corollary 4.6.
Let $F$ and $Q$ be as in (4.6) and (4.7) and $H\subset \R^4$ be the $F$ invariant subspace corresponding to the irreducible factor $x^3-x^2-2x -1$ in the characteristic polynomial for $F$.  Then both $\{{\bf s}(C^+_0), {\bf s}(C^+_1), {\bf s}(f^*C^+_0)\}$ and $\{{\bf u}(C^-_0), {\bf u}(C^-_1), {\bf u}(f^*C^-_0)\}$ are bases for $H$.  The action of $F$ on $H$ relative to either basis is given by the matrix (2.7).  If ${\bf s} = a_1{\bf s}(C^+_0)+a_2{\bf s}(C^+_1)+a_3{\bf s}(f^*C^+_0)$ and ${\bf u} = b_1{\bf u}(C^-_0) + b_2{\bf u}(C^-_1)+b_3{\bf u}(f^*C^-_0)$ are vectors in $\R^4$, then
${\bf s}^T Q{\bf u}$ is given by the right side of equation of equation (2.8).

This corollary motivates the following definition.  We let $\cC^+_R$ be the set of complex algebraic curves $\Gamma\subset {\bf C}^2$ satisfying the conditions:
\item{(i)} $\Gamma$ is linearly equivalent to $n_0 C^+_0+n_1 C^+_1 + n_2f^*( C^+_0)$ for some $n_0,n_1,n_2\in{\bf N}$.
\item{(ii)} ${\bf s}(\Gamma\cap{\bf R}^2) \ge n_0{\bf s}(C_0^+) + n_1{\bf s}(C_1^+) + n_2{\bf s}(f^*C_0^+)$.

\noindent We define $\cC^-_R$ in an analogous fashion, using the basis $\{ C^-_0, C^-_1,f_* C^-_0\}$.

\proclaim Proposition 4.7.  $\cC^\pm_R\ne\emptyset$,  $f^*\cC^+_R\subset \cC^+_R$, and $f_*\cC^-_R\subset\cC^-_R$.

\noindent{\it Proof.}  For $t>0$ we set $\Gamma=\{y=t\}$.  Then $\Gamma\sim C^{+}_{0}$ and $\Gamma\cap{\bf R}^2$ contains a 1s-arc, so we see that $\Gamma\in\cC^+_R$.  Hence $\cC^+_R$ and, similarly, $\cC^-_\R$ are non-empty.  The remaining conclusions follow immediately from the definition of $\cC^\pm_R$, equation (2.7) and Corollary 4.6.
\bigskip

\proclaim Corollary 4.8.  If $W\in\cC^+_R$, $W'\in\cC^-_R$, then
\item{(i)}  All intersection points of $W\cap W'$ are real.
\item{(ii)}  Each s-subarc of $W\cap{\bf R}^2$ intersects each u-subarc of $W'\cap {\bf R}^2$ transversally in exactly $Q_{i,j}$ points, and all points of $W\cap W'$ arise in this fashion.

\noindent{\it Proof.} The argument here is the same as the one used in the second half of the proof of Proposition 4.5.
Namely, by the definition of $\cC^\pm_R$, equation (2.8), and Corollary 4.6, we have that the complex intersection
number $W\cdot W'$ is bounded above by the number of distinct real intersections $W\cap W'\cap \cX_R$.  The only
way this can happen is if the two numbers are equal.  That is, all intersections between $W$ and $W'$ are the real
intersections between u-arcs and s-arcs given by Proposition 4.3, and no further intersections are possible.
\bigskip

\noindent{\bf \S5.  Coding Orbits. }  Here we establish a coding for $f$-orbits.  We will consider orbits which lie in the union of the open rectangles $R^+_j$ for $j=1,2,3,4$.  We say that a transition $(j,k)$ is {\it admissible} if  $fR_j^+\cap R_k^+\ne\emptyset$.  From Figure 4.3 and Equation (4.5), we see that the admissible transitions are given by the matrix $F$ in equation (4.6); that is, $(j,k)$ is admissible if and only if $F_{j,k}=1$.  We say that a word $w=w_{-n}\cdots w_{m}$ on the letters $1,2,3,4$  is {\it admissible} if  $F_{s_j,s_{j+1}}=1$ for all $j$.   We say that a sequence $p_{-n}\cdots p_{m}$ of points in $\cX_R$ is an {\it orbit segment} if $p_j\notin\cI(f)$ and $fp_j=p_{j+1}$ for $-n\le j<m$.  We say that this orbit segment is coded by the word $w$ if $p_{j}\in R_{w_j}^{+}$ for $-n\le j\le m$.

\proclaim Proposition 5.1.  Let $w=w_{-n}\cdots w_{m}$ be a finite word.  If $w$ is not admissible, then there is no orbit segment coded by $w$.  On the other hand, if $w$ is admissible, and $\gamma$ is an s-arc in $R_{w_m}$, then $f^{-m-n}\gamma$ contains an s-arc $\gamma'$ such that $w$ codes $p$, $fp$, \dots, $f^{m+n}p$ for every $p\in\gamma'$ except the endpoints.  In particular, $w$ defines a surjection from the finite orbit segments to the finite admissible sequences.

\noindent{\it Proof. }  If there is an orbit segment $\dots,p_{j},p_{j+1},\dots$ that is coded by $w$, then $p_{j+1}\in R^+_{w_{j+1}}\cap f(R^+_{w_{j}})$, so $w$ must be admissible.

Now suppose that $w$ is admissible.  An s-arc $\gamma$ in $R^{+}_{w_{m}}$ is of type $w_{m}s$.  By Proposition 4.1,  $f^{-1}$ contains a $j$s-arc $\gamma'$ for each $j$ such that $F_{j,w_{m}}=1$.  We choose a $w_{m-1}$s-subarc $\gamma_{1}$ of $f^{-1}\gamma$.  By the same argument, we can choose a $w_{m-2}$s-subarc $\gamma_{2}$ of $f^{-1}(\gamma_{1})$.  We continue in this way until we reach a $w_{-n}$s-arc $\gamma_{n+m}\subset f^{-n-m}\gamma_{n+m-1}$.  This arc satisfies the second conclusion of the Proposition.
\medskip

We let $\Sigma_F\subset\{1,2,3,4\}^{\bf Z}$ denote the space of bi-infinite sequences $w$ which are admissible.  We let $\sigma$ denote the shift operator on $\Sigma_F$, where $\sigma(w)=z$ means that $z_n=w_{n-1}$ for all $n\in{\bf Z}$.  The sequences $\overline{3}$ and $\overline{4}$ are the fixed points of $\sigma$, and the 2-cycles are $\sigma(\overline{12})=\overline{21}$ and $\sigma(\overline{34})=\overline{43}$.  The space $\Sigma_F$ inherits the infinite product topology from $\{1,2,3,4\}^{\bf Z}$, and thus it is compact and totally disconnected.  The stable manifold of the 2-cycle $\{\overline{12},\overline{21}\}$ consists of the set of all sequences $*\overline{12}$, which are ultimately repeating `$12$' to the right.  The sequences which are ultimately repeating `$12$' to the left are the unstable manifold of this 2-cycle.  Let us use the notation
$$W^{s/u}({\bf 2}):=W^s(\overline{12},\overline{21})\cup W^s(\overline{34},\overline{43})\cup W^u(\overline{12},\overline{21})\cup W^u(\overline{34},\overline{43})$$
for the union of the stable and unstable manifolds of the 2-cycles in $\Sigma_F$.

Let us define $\Omega=\bigcap_{n\in{\bf Z}}f^n(R^+_1\cup R^+_2\cup R^+_3\cup R^+_4)$.  By the construction of the sets $R^+_j$, we see that if $p\in\Omega$, then $f$ is a local diffeomorphism at $f^np$ for all $n\in{\bf Z}$.  Thus $\Omega$ is a natural domain of definition for the ``coding map'' that takes the orbit of a point $p\in\Omega$ and assigns its itinerary:
$$w:\Omega\to\Sigma_F, \ \ \Omega\ni p\mapsto (w_n)_{n\in{\bf Z}}, \ \ f^np\in R^+_{w_n}, \ \forall n\in{\bf Z}. $$
The coding map is a semi-conjugacy from $\Omega$ to $\Sigma_F$.

\proclaim Theorem 5.2.  The map $w:\Omega\to\Sigma_F$ is continuous, and $w(\Omega)\supset\Sigma_F-W^{s/u}({\bf 2})$.

\noindent{\it Proof. }  Continuity of the coding map follows from the facts that $\Omega\subset R^{+}_{1}\cup R^{+}_{2}\cup R^{+}_{3}\cup R^{+}_{4}$ and that $f$ is uniformly continuous on compact subsets of each $R^{+}_{j}$.

Now let $s\in\Sigma_F$ be given, and define
$$M_{k}=\{p\in\cX_{R}:f^{j}p\in R^{+}_{s_{j}}{\rm\ for\ }-k\le j\le k\}=\bigcap_{j=-k}^{k}f^{-j}R^{+}_{s_{j}}.$$
Clearly $M_{k}\supset M_{k+1}$, and by Proposition 5.1, $M_{k}\ne\emptyset$ for each $k$.  Thus $M:=\bigcap_{k\ge0}\bar M_{k}$ is a nonempty, compact set.

If $M\subset\Omega$, then $w(M)=s$.  Suppose instead that there is a point $p\in M-\Omega$.  We will complete the
proof by showing that $s\in W^{s/u}({\bf 2}).$

\proclaim Lemma 5.3.  If $p\in M-(\Omega\cup V_0\cup V_2)$, then one or both of the following is true:
\item{} $p\notin \cI(f^n)$ for any $n>0$, and $f^n(p)\in \bigcup E_j^s$ for $n>0$ large enough;
\item{} $p\notin \cI(f^{-n})$ for any $n>0$, and $f^{-n}(p)\in\bigcup E^u_j$ for $n>0$ large enough.

\noindent{\it Proof.}
Suppose first that $p\in V_1\cup V_3\cup V_4\cup V_5$.  Theorem 1.1 tells us that $p$ cannot be indeterminate for both a forward and a backward iterate of $f$.  If, for example, $p\notin\bigcup_{n>0} \cI(f^n)$, then $f^n(p)$ will be well-defined and contained in $\bigcup E^s_j$ for $n>0$ large enough.  Likewise, $p\notin\bigcup_{n>0} \cI(f^{-n})$ implies that $f^{-n}(p)\in \bigcup E^u_j$ for $n$ large.

Now suppose that $p\in\cX_R-\bigcup V_j$.  Then $p\notin \cI(f^n)$ for any $n\in\Z$.  Since $p\notin\Omega$, we can choose $k\in\Z$, with $|k|>0$ minimal so that $f^k(p)\notin R^+_1\cup\dots R^+_4$.  Without loss of generality, say $k>0$. Since $p\in \overline{M_k}$, it
follows that $f^k(p)\in \bigcup_{j=1}^4 bR^+_j - \bigcup V_j \subset \cC(f)$.  Hence $f^n(p)\in \bigcup E^s_j$ for all $n>k$. \medskip

Returning to the proof of Theorem 5.2 with $p\in M-\Omega$, we have by definition of $M$ that any finite piece $f^{-k}(p),\dots, f^k(p)$ of the orbit of $p$ must be approximated by $f^{-k}(q)\in R^+_{s_{-k}},\dots, f^k(q)\in R^+_{s_k}$ for some $q\in M_k$.  Suppose that $p$ satisfies the first alternative in Lemma 5.3: $f^n(p)\in \bigcup E^s_j$ for $n>0$ large.  More specifically, by (1.2) either $f^n(p)$ alternates between $E^s_4$ and $E^s_5$, or $f^n(p)$ alternates between $E^s_1$ and $E^s_3$.  Take the first case, for example.  The orbit of $q$ is admissable and $f(R_1^+)\cap R_1^+ = R_2^-\cap R_1^+ = \emptyset$.  Hence it is apparent from Figures 4.1 and 4.3 that the only way $f^n(q)$ can remain near $f^n(p)$
is if $f^n(q)\in R_1^+$ when $f^n(p)\in E^s_4$, and $f^n(q)\in R_2^+$ when $f^n(p)\in E_5^s$.  That is,
$s=w(M)\in W^s(\overline{12},\overline{21})$.  Likewise, if $f^n(p)$ alternates between $E^s_1$ and $E^s_3$, similar
analysis shows that $s\in
W^s(\overline{34},\overline{43})$.

If $p$ satisfies the second alternative in Lemma 5.3, then we can repeat our arguments to conclude that
$s$ belongs to $W^u(\overline{12},\overline{21})$ or $W^u(\overline{34},\overline{43})$.

Finally, if $p\in V_0\cup V_2$, then the orbit of $p$ is contained in the `vertices' of the hexagon in Figure 4.3. Combining the fact that the orbit of $p$ is shadowed by an admissable orbit with equations (1.2) and (4.5), one deduces that in fact there are only two possibilities to consider:  $f^n(p)$ lies in the piece of $V_0$ corresponding to the left/right vertex of the hexagon for all $n\in\Z$, or $f^n(p)$ lies in the portion of $V_2$ corresponding to the point(s) $V_1\cap V_3$ in the hexagon.  In the first case, we can repeat the arguments above to conclude that $s$ is $\overline{12}$ or $\overline{21}$.  In the second case, we obtain that $s$ is $\overline{34}$ or $\overline{43}$.

\bigskip\noindent{\bf\S6.  Orbits Attracted to Infinity. }
By  ``infinity'' we will mean the set
$$
\cX_{R}-{\bf R}^{2}= \cX_R\cap(V_{0}\cup V_{2}\cup V_{3}\cup V_{4}\cup V_{5}).
$$
Recall from (1.2) that $f^{2}$ acts as a translation on $V_{j}$ for $j=3,4,5$.  The curves $V_{0}$ and $V_{2}$ are fixed pointwise by $f^{2}$.  In fact, each point of $V_{0}\cup V_{2}$ is a parabolic fixed point for $f^{2}$.  In this section we will show (Theorem 6.3) that each point $p\in{\bf R}^{2}-\Omega$ is attracted to infinity in forward time or backward time, or both.

\proclaim Lemma 6.1.  $f(R^{-}_{6}\cup R^{-}_{7})\subset R^{+}_{5}$, and $f(R^{+}_{5})\subset R^{-}_{7}$.

\noindent{\it Proof. }  From Figure 4.3, we see that $R^{-}_{6}\cup R^{-}_{7}\subset R^{+}_{1}$.  So (4.5) tells us that
$f(R^-_6\cup R^-_7)\subset R^-_2$.  Further, the finite part of the boundary of $\overline{f(R^-_6\cup R^-_7)}$ consists
of an arc $\gamma\subset f(C_0^-)$, and since it joins $E^s_3$ to $E^s_4$, $\gamma$ is not an s-arc.
Corollary 4.8 therefore allows us to conclude that $\gamma\cap C_1^+ = \emptyset$.  This proves that
$f(R^{-}_{6}\cup R^{-}_{7})\subset R^{+}_{5}$.

The second assertion is immediate from (4.5).
\medskip
\epsfysize=1.9in
\centerline{ \epsfbox{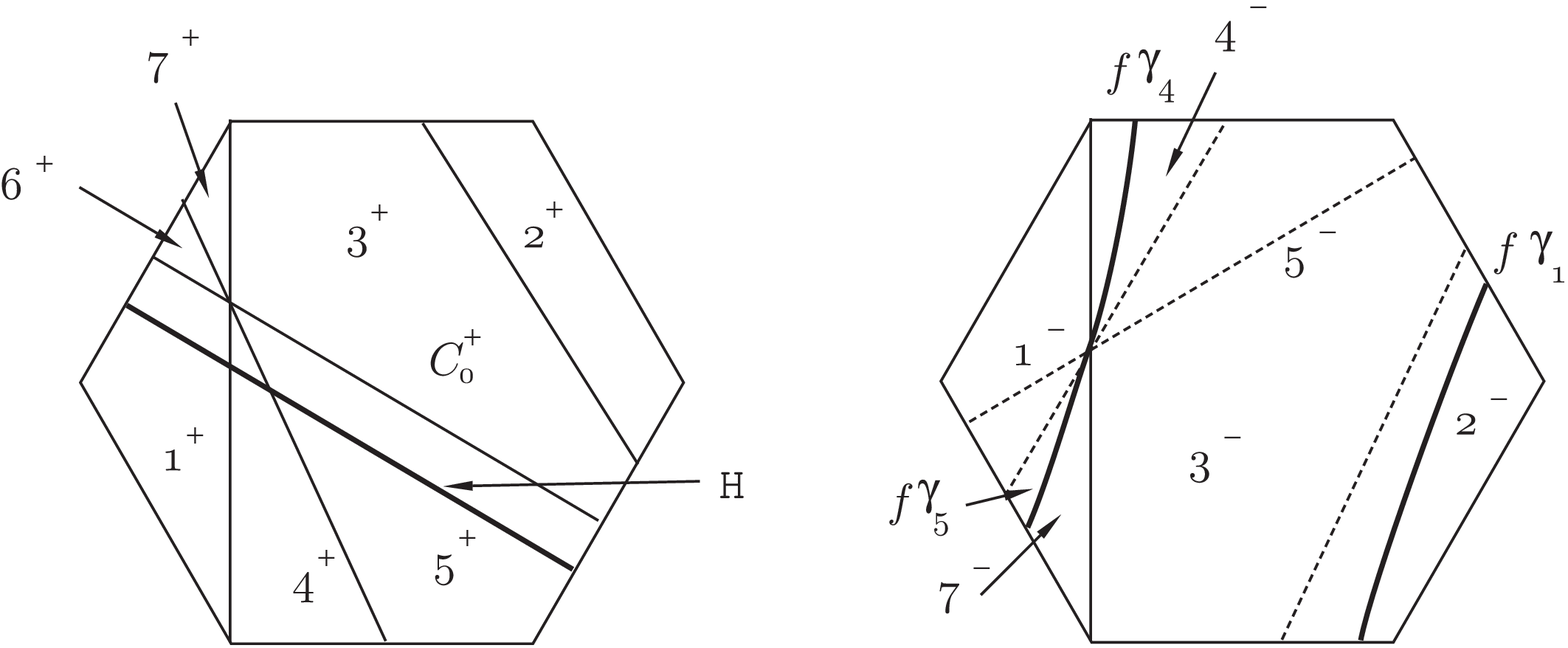}  }
\centerline{\rm Figure 6.1.  Forward image of a horizontal line.}
\medskip

Lemma 6.1 implies that $f^2$ maps $R^+_5$ into itself.  The next results shows that the action $f^2|R^+_5$ moves points down
by a definite amount when $f^2$.  We let $\phi(x,y)=y$ be projection onto the vertical coordinate.

\proclaim Lemma 6.2.  If $p\in R^+_5$, then $\phi(f^{2}p)-\phi(p)\le 2a$.

\noindent {\it Proof. }   Given $p\in R^+_5$, let $H$ be the horizontal line $\{(x,y):y=\phi(p)\}$.
For $a<a'<0$, let $H'=\{(x,y):\phi(x,y)=\phi(p)+2a'\}$.  Then as complex curves, $H$ and $H'$ are linearly equivalent to $C_0^+$.  So from (2.8), we compute  that $H\cdot f^2_* H' = 3.$
\medskip
\epsfysize=1.7in
\centerline{ \epsfbox{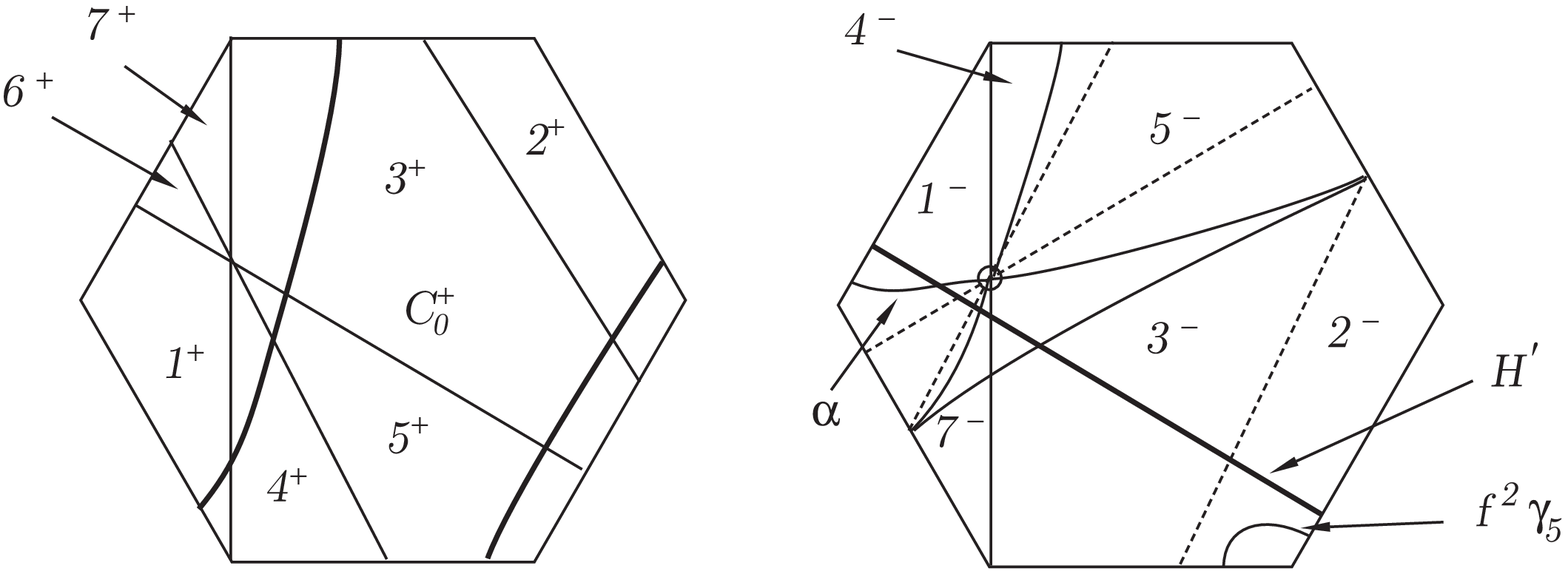}  }
\centerline{\rm Figure 6.2.  Placement of $\alpha$ and $f^2\gamma_5$.}
\medskip
Since $p\in R^+_5$, we have $\phi(p)<0$, and thus $H$ crosses the regions $R^+_j$ as shown on the left side of Figure 6.1.  In particular, $H$ consists of arcs $\gamma_1=H\cap R^+_1$ and $\gamma_4=H\cap R^+_4$, followed by $\gamma_5=H\cap R^+_5$.  Applying (4.5), we see that $f\gamma$ is as in the right side of Figure 6.1.  We repeat this process with $f(\gamma)$
(left side of Figure 6.2) in place of $\gamma$ to obtain $f^2(\gamma)$ (right side of Figure 6.2).  As before, the results are determined by (4.5) and by the action of $f$ on $\cC(f)\cup V_0\cup \dots V_5$.  From the right side of Figure 6.2, it is clear that $H'$ must cross $f^2(\gamma)$ at least three times.  In particular, since there are only three (complex)
intersections, $H'$ cannot intersect the small hook $\gamma' = f^2(\gamma)\cap R^-_2$ near the lower right corner of the hexagon.  It follows that $\gamma'$ lies completely below $H'$.  On the other hand, the facts that $f^2(R_5^+)\subset R_5^+$ and $p\in\gamma\cap R_5^+$ imply that $f^2(p)\in \gamma'$.  We conclude that $\phi(f^2(p)) < \phi(p)+2a'$, thereby finishing
the proof.

\proclaim Theorem 6.3.   If $p\in\cX_R-\Omega$, then  $f^np$ tends to infinity as either $n\to+\infty$, or $n\to-\infty$, or both.

\noindent{\it Proof. }  Since $p\notin\Omega$, some forward or backward iterate $f^n(p)$ leaves $\bigcup_{j=1}^4R^+_j$.  Suppose first that $n\geq 0$.  If $f^n(p)\in V_0\cup\dots\cup V_5\cup \cC(f)$, then the
assertion is immediate from (1.2).  Otherwise, $f^n(p)\in R_5^+$, because the analogue of Lemma 6.1 prevents
$f^n(p)\in R_6^+\cup R_7^+$.  In this case, $f^{n+k}(p)$ alternates between $R_5^+$ and $R_6^-\cup R_7^-$ for all
$k\geq n$, and Lemma 6.2 tells us that the forward orbit tends to infinity.  When $n<0$, then the same argument
shows that either $f^n(p)\in V_0\cup\dots\cup V_5\cup \cC(f^{-1})$ or $f^n(p)\in R_6^+\cup R_7^+$, and in both
cases the backward orbit of $p$ tends to infinity.

\proclaim Theorem 6.4.  Let $a<-1$, $b\ne0$, and let $p,q\in{\bf R}^{2}$ be saddle points.  Let $W^{s}_{C}(p)\subset\cX$ denote the complex stable manifold through $p$, and let $W^{u}_{C}(q)$ denote the complex unstable manifold through $q$.  Then $W^{s}_{C}(p)\cap W^{u}_{C}(q)\subset{\bf R}^{2}$.

\noindent{\it Proof.}  We let $D^u$ denote a disk inside $W^u(p)$ which contains $p$, and we let $D^s\ni q$ denote a disk inside $W^s(q)$.  Note that neither disk is contained in any algebraic curve, since such a curve would automatically be
invariant under some iterate of $f$, and the only invariant curves are $V_0, \dots ,V_5$.  We seek to approximate
$D^u$ and $D^s$ by curves in $\cC^-_R$ and $\cC^+_R$, respectively.

By Theorem 6.3 we have $p,q\in\bigcup_{j=1}^4R^+_j$.  Further, if $p\in R^+_1$ then $fp\notin R^+_1$, so we may assume that $p,q\notin R_1^+$.  Note that curves of the form $\{(x,y)\in {\bf C}^2:y=s\}$, $s>0$, each belong to $\cC^+_\R$ and
collectively fill the regions $R_2^+\cup R_3^+$.  Similarly, curves of the form $H = \{(x,y)\in{\bf C}^2:(x-1)(y-s)=s\}$,
$s > 1$, belong to $\cC^+_\R$ and fill $R_4^+$.  Therefore, we may choose some curve $H\in\cC^+_\R$ that intersects $D^u$ transversally in a point $p'$.  Likewise, we may choose a curve $H'\in \cC^-_R$ that intersects $D^s$ transversally at a point $q'\in D^s$.

Now suppose that $r\in W^s_C\cap W^u_C$.  Then $r\in f^{-n}(D^s)\cap f^n(D^u)$ for $n$ large.  Shrinking $D^s$ and $D^u$
if necessary, we can apply the Lambda Lemma (see [R, \S5.11]) to obtain disks $D^s_0\subset H$ and $D^u_0\subset H'$
such that $f^n D^u_0$ is $C^1$ close to $f^n(D^u)$ and $f^{-n}D^s_0$ is $C^1$ close to $f^{-n}(D^s)$.  We conclude,
then, that $f^n D^u_0$ intersects $f^{-n}D^s_0$ in a point near $r$.  On the other hand, $f^nD^u_0$ is contained in an element of $\cC^-_R$, and $f^{-n}D^s_0$ is contained in an element of $\cC^+_R$.  It follows from Corollary 4.6 that their
intersection must be real.  Thus $r\in{\bf R}^2$.

\bigskip\noindent{\bf \S7.  Invariant Measure. }  Let us recall some results about ergodic theory for birational maps on complex surfaces.  We suppose that $f:X\to X$ is a birational map satisfying the conditions in Theorem 1.1.  We suppose, in addition, that the spectral radius of $f^*$ on ${\rm Pic}(X)$ is $\rho>1$.   By [DF] there are positive, closed (1,1)-currents $T^\pm$ with the property that $f^*T^+=\rho T^+$ and $f_*T^-=\rho T^-$.  The paper [BD2] introduced a quantitative hypothesis on the speed of approach of $f^{-n}\cI(f)$ to $\cI(f^{-1})$ and showed that when it holds the product $\mu := T^+\wedge T^-$ is a well-defined ergodic invariant measure.  By [D], this measure has entropy $\log\rho$.  If (4.4) holds, then our maps $f_{a,b}$ satisfy the hypothesis given in [BD], for in this case $f^{-n}\cI(f)$ remains at a bounded distance from $\cI(f^{-1})$.

For $t\in{\bf R}$ we set $H_{t}=\{(x,y)\in{\bf C}^{2}:y=t\}$, and we let $[H_{t}]$ denote the current of integration over $H_{t}$.  For $0<t_{1}<t_{2}$ we define the currents
$$
\omega^{+} = (t_{2}-t_{1})^{-1}\int_{t_{1}}^{t_{2}} [H_{t}]\,dt
$$
and $\omega^{-}=\tau^{*}\omega^{+}$.  The normalization is chosen so that $\omega^+$ and $\omega^-$ represent the fundamental classes of $H_t$ and $\tau^*H_t$, respectively.  By the Poincar\'e-Lelong formula, the potential of $[H_{t}]$ is $(2\pi)^{-1}\log|y-t|$, and so the potential of $\omega^{+}$ is
$$(t_{2}-t_{1})^{-1}\int_{t_{1}}^{t_{2}}\log|y-t|\,dt,$$
which is continuous.  The pull-back of $[H_{t}]$ is $f^{n*}[H_{t}]=[f^{-n}H_{t}]$, and the pull-back of $\omega^{+}$ is given by
$$
f^{n*}\omega^{+}=(t_{2}-t_{1})^{-1}\int_{t_{1}}^{t_{2}}[f^{-n}H_{t}]\,dt.
$$
The wedge product of varieties is the sum of point masses over the intersection:
$$
[f^{-n}H_{s}]\wedge[f^{n}\tau H_{t}] = [f^{-n}H_{s}\cap f^{n}\tau H_{t}] = \sum_{p\in f^{-n}H_{s}\cap f^{n}\tau H_{t}} \delta_{p}.
$$
We may integrate this to obtain
$$
f^{n*}\omega^{+}\wedge f^{n}_{*}\omega^{-}= (t_2-t_1)^{-2}\int_{t_{1}}^{t_{2}} ds \int_{t_{1}}^{t_{2}} dt \, [f^{-n}H_{s}\cap f^{n}\tau H_{t}]. \eqno(7.1)
$$
Now we apply [BD2]:  since the currents $\omega^{\pm}$ have continuous potentials, the invariant measure $\mu$ is a positive
multiple of:
$$
\lim_{n\to\infty} \rho^{-2n}f^{n*}\omega^{+}\wedge f^{n}_{*}\omega^{-}.\eqno(7.2)
$$

\proclaim Theorem 7.1.  If $a<-1$ and $b\ne0$, then $\mu$ is carried by $\Omega\subset{\bf R}^2$.

\noindent{\it Proof. }  As in the proof of Proposition 4.7, we have $H_{s}\in\cC^{+}_{R}$ and $\tau H_{t}\in\cC^{-}_{R}$ for $s,t>0$.  Thus by Corollary 4.8, $f^{-n}H_{s}\cap f^{n}\tau H_{t}\subset{\bf R}^{2}$.  In fact since the intersections come from crossings between s-arcs and u-arcs, they all lie in $R_1^+\cup\dots\cup R_4^+$.  Thus the measures in (7.1) are supported in $\bigcup_{j=1}^4 R_j^+$.  By (7.2) and the fact (see [BD]) that $\mu$ does not charge algebraic curves, we conclude that $\mu$ is supported on $\bigcup_{j=1}^4 R_j^+ \subset \cX_{R}$. Since $\mu$ is invariant, $\mu$ must be carried by $\Omega$.
\bigskip

Let us recall the construction of the unique measure $\nu$ of maximal entropy on $\Sigma_F$ (cf.\ [LM] for further details). For an admissible word $w'=w'_{-N}\cdots w'_{N}$ we define the cylinder set
$$
C(w'):=\{w\in\Sigma_F:w_{j}=w'_{j}, -N\le j\le N\}.
$$
The cylinder sets are open (and closed) and generate the topology of $\Sigma_F$.  We define
$$
\nu(C(w')):=\lim_{n\to\infty}{\#\{w_{-n}\cdots w_{n}{\rm\ admissible }:w_{j}=w'_{j}, -N\le j\le N\} \over
\#\{w_{-n}\cdots w_{n}{\rm\ admissible }\} }.\eqno(7.3)
$$
The system $\Sigma_F$ is transitive since all entries of $F^{3}$ are strictly positive.  Thus the limit in (7.3) does not
change if we restrict on the right side to words $w$ with prescribed first and last symbols.  We will take advantage
of this flexibility below.

\proclaim Theorem 7.2.  If $a<-1$ and $b\ne0$, then the map $w:\Omega\to\Sigma_F$ sends $\mu$ to $\nu$, i.e., $w_*\mu=\nu$.

\noindent{\it Proof. }  Let $w'=w'_{-N}\cdots w'_{N}$ be an admissible word.  We will show that $\mu(w^{-1}C(w'))=\nu(C(w'))$.  Define
$$
R(w'):=\bigcap_{j=-N}^{N} f^{-j}R^{+}_{w'_{j}}.
$$
Then $R(w')\cap\Omega$ is open and closed in $\Omega$ because $R^+_{w'_j}\cap\Omega$ is open and closed for each $j$. By Theorem 7.1
$$
\mu R(w')  = \mu(R(w')\cap\Omega) = \mu(w^{-1}C(w')).
$$
So by (7.2) it will suffice to show that
$$
\lim_{n\to\infty}\mu_{n}R(w')=\nu(C(w')).\eqno(7.4)
$$
where $\mu_{n} = c_n f^{n*}\omega^+\wedge f^n_*\omega^-$, and $c_n>0$ is chosen so that $\mu_n$ has unit mass.

Fix numbers $0<s<1$ and $t>0$, and let $H_{s}=\{y=s\}$ and $\tau H_{t}$ be individual lines appearing in the integrals that
define $\omega^+$ and $\omega^-$.  Let $S_n = f^{-n}(H_s) \cap f^n(\tau H_t)$ and
$\tilde S_n$ be the set of admissible words $w_{-n}\dots w_n$ for which $w_n = 3$ and $w_{-n} = 3$ or $4$.

\proclaim Lemma 7.3.
The coding map sends $S_n$ bijectively onto $\tilde S_n$.

\noindent{\it Proof. }
As Figures 4.3 and 4.4 make clear, the restriction on $s$ implies that $H_{s}$ contains a 3s-arc in $R^{+}_{3}$.
Similarly, $\tau H_t$ contains a 3u-arc.  Proposition 5.1 therefore implies for each admissible word
$w=w_0\cdots w_{2n-1}3$ ending with 3, that
$f^{-2n}(H_s)$ contains at least one s-arc coded by $w$.  Among these, the 3s-arcs and 4s-arcs each intersect
the 3u-arc in $\tau H_t$ at least once.  This gives a total of at least $F^{2n}_{3,3}+F^{2n}_{4,3}$ points in
$f^{-2n}(H_s)\cap \tau H_t$, each coded by a distinct admissible word.  Applying $f^n$ to points in
$f^{-2n} H_s\cap \tau H_t$ and shifting the codings $n$ symbols to the left, we arrive at
a subset of $S_n$ mapped bijectively onto $\tilde S_n$ by coding.

To see that there are no other points in $S_n$, we count complex intersections.  Since $H_s\in\cC^+_R$ and $H_t\in \cC^-_R$
are linearly equivalent to $C_0^+$ and $C_0^-$, respectively, we use Propositions~4.3 and 4.5 and Corollary 4.8 to compute
$$
\eqalign{&f^{n*}H_s\cdot f^n_* \tau H_t  =  H_s \cdot f^{2n}_*\tau H_t
= \cr
=&(0,0,1,0) \cdot QF^{2n}\pmatrix{ 0\cr 0\cr 1\cr 0\cr} = (0,0,1,1)\cdot  F^{2n}\pmatrix{ 0\cr 0\cr 1\cr 0\cr} = F^{2n}_{3,3} + F^{2n}_{4,3}.\cr}
$$
The number of real intersections is bounded above by the number of complex intersections, and the computation
shows therefore that the two quantities are equal.
\medskip

Returning to the proof of Theorem 7.2, we have from Lemma 7.3 that the number of elements in $S_n$ is independent
of $s$ and $t$.  Thus from (7.1) and Lemma 7.3
$$
\mu_n R(w') = {\# S_n\cap R(w')\over \# S_n} = {\# \tilde S_n\cap C(w') \over \# \tilde S_n}.
$$
In light of (7.3) and the ensuing remarks, we conclude that (7.4) holds.
\proclaim Corollary 7.4.  The entropy of $f$ is $\log\rho$.

\noindent  {\it Proof. }   By [DS] the entropy of $f$ is bounded above by $\log\rho$.  The other inequality is seen because by Theorem 7.2, the entropy of $f$ is at least as great as the entropy of $\nu$, which is equal to $\log\rho$. 

\bigskip\centerline{\bf References}

\item{[A1]}    N. Abarenkova,  J-C. Angl\`es d'Auriac, S. Boukraa,
    S. Hassani and J-M. Maillard,  From Yang-Baxter equations to dynamical zeta functions for birational transformations.  Statistical physics on the eve of the 21st century, 436--490, Ser. Adv. Statist. Mech., 14, World Sci. Publishing, River Edge, NJ, 1999.

\item{[A2]}  N. Abarenkova,  J-C. Angl\`es d'Auriac, S. Boukraa,
    S. Hassani and J-M. Maillard,
    Rational dynamical zeta functions for
     birational transformations.
    Physica {\bf A 264} (1999) pp. 264--293.  chao-dyn/9807014.

\item{[A3]} N. Abarenkova, J-C. Angl\`es d'Auriac,
   S.~Boukraa, S. Hassani and  J-M. Maillard,
   Topological entropy and complexity for discrete dynamical
   systems.     Phys. Lett. {\bf A 262} (1999) pp. 44--49.  chao-dyn/9806026 .

\item{[A4]}  N. Abarenkova, J-C. Angl\`es d'Auriac, S. Boukraa and J-M. Maillard, Growth complexity spectrum of some discrete dynamical systems.  Phys. D 130 (1999), no. 1-2, 27--42.

\item{[A5]}  N. Abarenkova, J-C. Angl\`es d'Auriac, S. Boukraa and J-M. Maillard, Real topological entropy versus metric entropy for birational measure-preserving transformations.  Phys. D 144 (2000), no. 3-4, 387--433.

\item{[A6]}  N. Abarenkova, J-C. Angl\`es d'Auriac, S. Boukraa, S. Hassani and J-M. Maillard, Real Arnold complexity versus real topological entropy for birational transformations.  J. Phys. A 33 (2000), no. 8, 1465--1501.

\item{[A7]}  N. Abarenkova, J-C. Angl\`es d'Auriac, S. Boukraa, S. Hassani and J-M. Maillard, Topological entropy and Arnold complexity for two-dimensional mappings.  Phys. Lett. A 262 (1999), no. 1, 44--49.

\item{[BD1]} E. Bedford and J. Diller,  Real and complex dynamics of a family of birational maps of the plane: the golden mean subshift, American J. of Math., to appear.

\item{[BD2]} E. Bedford and J. Diller,  Energy and invariant measures for birational surface maps.  Duke Math.\ J., to appear.

\item{[BTV]} M. Bernardo, T. Truong and G. Rollet, The discrete Painlevé I equations: transcendental integrability and asymptotic solutions. J. Phys. A 34 (2001), no. 15, 3215--3252.

\item{[BV]} M. Bellon and C. Viallet, Algebraic entropy. Commun. Math. Phys. 204, 425--437 (1999)

\item{[BHM1]}  S. Boukraa,   S. Hassani and  J-M. Maillard,  New integrable cases of a Cremona transformation: a finite-order orbits analysis, Physica A 240 (1997), 586--621.

\item{[BHM2]}  S. Boukraa,   S. Hassani and  J-M. Maillard,  Product of involutions and fixed points, Alg.\ Rev.\ Nucl.\ Sci., Vol. 2 (1998), 1--16.

\item{[BMR]}  S. Boukraa, J-M. Maillard, and G. Rollet, Almost integrable mappings,  Int.\ J. Mod. Phys. {\bf B8} (1994), pp. 137--174.

\item{[DF]}  J. Diller and C. Favre,  Dynamics of bimeromorphic maps of surfaces, American J. of Math., 123 (2001), 1135--1169.

\item{[DS]} T-C. Dinh and N. Sibony,  Une borne sup\'erieure pour l'entropie topologique d'une application rationnelle.  preprint.

\item{[D]}  R. Dujardin,  Laminar currents and birational dynamics.  preprint

\item{[FS]} J-E. Forn\ae ss and N. Sibony, Complex dynamics in higher dimension, II, {\sl Modern Methods in Complex Analysis},  Ann.\ of Math.\ Studies, vol. 137, Princeton U. Press, 1995, p. 135--182.

\item{[GH]} P. Griffiths and J. Harris, Principles of algebraic geometry. Wiley Classics Library, 1994.

\item{[LM]}  D. Lind and B. Marcus,  An introduction to symbolic dynamics and coding. Cambridge University Press, Cambridge, 1995.

\item{[R]} C. Robinson,  Dynamical systems. Stability, symbolic dynamics, and chaos. Second edition. Studies in Advanced Mathematics. CRC Press, Boca Raton, FL, 1999.

\bigskip
\rightline{Indiana University}

\rightline{Bloomington, IN 47405}

\rightline{\tt bedford@indiana.edu}

\medskip
\rightline{University of Notre Dame}

\rightline{Notre Dame, IN 46556}

\rightline{\tt diller.1@nd.edu}

\bye